\documentclass[12pt]{article}
\usepackage{mathrsfs,amssymb,amsmath,amsthm,tikz,multirow,graphicx}
\usetikzlibrary{calc,arrows}
\usepackage{caption}
\usepackage{graphicx}
\usepackage{float} 
\usepackage{subcaption}

\usepackage[hidelinks, hyperindex]{hyperref}

\hyperref[sec:function]{}

\newcommand{\arcThroughThreePoints}[4][]{
\coordinate (middle1) at ($(#2)!.5!(#3)$);
\coordinate (middle2) at ($(#3)!.5!(#4)$);
\coordinate (aux1) at ($(middle1)!1!90:(#3)$);
\coordinate (aux2) at ($(middle2)!1!90:(#4)$);
\coordinate (center) at ($(intersection of middle1--aux1 and middle2--aux2)$);
\draw[#1] 
 let \p1=($(#2)-(center)$),
      \p2=($(#4)-(center)$),
      \n0={veclen(\p1)},       % Radius
      \n1={atan2(\y1,\x1)}, % angles
      \n2={atan2(\y2,\x2)},
      \n3={\n2>\n1?\n2:\n2+360}
    in (#2) arc(\n1:\n3:\n0);
}

\newcommand*\circled[1]{\tikz[baseline=(char.base)]{
             \node[shape=circle,draw,inner sep=0.5pt] (char) {#1};}}
             
\newcommand*\squared[1]{\tikz[baseline=(char.base)]{
             \node[shape=rectangle,draw,inner sep=1 pt] (char) {#1};}}

\title{Moduli Spaces of Pentagonal Subdivision Tilings}
\author{Jinjin Liang, Erxiao Wang\thanks{Corresponding author (wang.eric@zjnu.edu.cn).  Research was supported by National Natural Science Foundation of China NSFC-RGC 12361161603 and Key Projects of Zhejiang Natural Science Foundation LZ22A010003.}, 
Zhejiang Normal University \\
Min Yan\thanks{Research was supported by NSFC-RGC Joint Research Scheme N-HKUST607/23 and Hong Kong RGC General Research Fund 16305920.}, 
Hong Kong University of Science and Technology}

\date{}

\newcommand{\bb}{\mathbb}

\newtheorem{theorem}{Theorem}
\newtheorem{lemma}[theorem]{Lemma}

\newtheorem*{theorem*}{Theorem}

\theoremstyle{definition}

\theoremstyle{remark}

\numberwithin{equation}{section}

\begin{document}

\maketitle

\begin{abstract}
Pentagonal subdivision gives three families of edge-to-edge tilings of the sphere by congruent pentagons. Each family forms a two dimensional moduli space. We describe these moduli spaces in detail.


\end{abstract}

Sommerville \cite{sommerville} started the classification of edge-to-edge tilings of the sphere by congruent triangles in 1924, and Ueno and Agaoka \cite{ua} completed the classification in 2002. In \cite{ay1, awy,cly1,gsy,cly2,lw,lwqx1,lwqx2,wy1,wy2}, we completely classified edge-to-edge tilings of the sphere by congruent polygons. The tilings are the Platonic type, the earth map type, and several quadrilateral tilings of the sporadic type. 

The tilings allow up to two free parameters. Quite a number of tilings allow one free parameter, and the ranges of the parameters are well understood. There are five tilings allowing two free parameters: The tetrahedron $P_4$, the first quadrilateral earth map tiling $E_{\square}1$, and the three pentagonal subdivisions $PP_4,PP_8,PP_{20}$ of the Platonic solids. The moduli space of the tiling $P_4$ is all the spherical triangles such that the sum of three angles is $2\pi$. The moduli space of the tiling $E_{\square}1$ is described in \cite[Figure 51]{cly2} and \cite[Figure 13]{lwqx1}. The moduli spaces of the pentagonal subdivisions $PP_4,PP_8,PP_{20}$ are more complicated, and are discussed in detail in this paper.

The pentagonal subdivision divides (more precisely, replaces) each triangular face $F$ of the regular tetrahedron $P_4$, octahedron $P_8$ or icosahedron $P_{20}$ into three congruent pentagons with edge lengths $a,a,b,b,c$. See the subdivision scheme in the first of Figure \ref{div_triangle}. Actually we do not keep the edges of $F$ straight. The second of Figure \ref{div_triangle} is the perspective picture of the stereographic projection of the pentagonal subdivision of a face $F$ of the tetrahedron, while Figure \ref{pst} gives real 3D pictures of three pentagonal subdivision tilings. The angles of the regular triangle $F$ are $\frac{2}{n}\pi$. Here and throughout the paper, $n=3,4,5$ is reserved for the tetrahedron, octahedron, icosahedron.

\begin{figure}[htp]
\centering
\begin{tikzpicture}[scale=1.5]

%% scheme

\begin{scope}[xshift=-3cm]

\foreach \a in {0,1,2}
{
\begin{scope}[rotate=120*\a]

\draw
	(0,0) -- (30:0.577);
	
\draw[line width=1.2] 
	(30:0.577) -- (0:1) -- (-30:0.577);

\draw[dashed]
	(30:0.577) -- (90:0.577);

\end{scope}
}

\filldraw[fill=gray!50] 
	(30:0.577) circle (0.05);

\node at (50:0.3) {\small $a$};
\node at (-110:0.3) {\small $a$};
\node at (18:0.85) {\small $b$};
\node at (-18:0.85) {\small $b$};
\node at (-55:0.58) {\small $c$};

\end{scope}

%% actual

\foreach \a in {0,1,2}
{
\begin{scope}[rotate=120*\a]

\draw
	(0:1) arc (30:90:{sqrt(3)});

\coordinate (M) at (60:{sqrt(3)-1});
\coordinate (X) at (0:1);

\pgfmathsetmacro{\th}{18} 

\pgfmathsetmacro{\ph}{-atan((1/3)*(cot(\th+30)))}
\coordinate (A1) at ({\ph}:{sqrt(2)});
\coordinate (B1) at ({\ph+180}:{sqrt(2)});

\pgfmathsetmacro{\ps}{-atan((1/3)*(cot(\th+150)))}
\coordinate (A2) at ({\ps}:{sqrt(2)});
\coordinate (B2) at ({\ps+180}:{sqrt(2)});

\coordinate (V) at (36.8:0.595);
\coordinate (W) at (-42:0.95);
\coordinate (W2) at (120-42:0.95);

\arcThroughThreePoints[line width=1.2]{X}{A1}{V};
\arcThroughThreePoints[line width=1.2]{W}{A2}{X};
\arcThroughThreePoints[dashed]{V}{M}{W2};

\draw 
	(0,0) -- (V);

\draw[gray]
	({-sqrt(3)+1},0) -- (1,0);

\end{scope}
}

%\fill (36.8:0.595) circle (0.05);

\fill 
	(1,0) circle (0.05);
\filldraw[fill=white]
	(0,0) circle (0.05);
\filldraw[fill=gray!50] 
	(36.8:0.595) circle (0.05);

\node at (0:1.15) {$B$};
\node at (120:1.15) {$B'$};
\node at (220:0.18) {$A$};
\node at (20:0.5) {$V$};
\node at (-65:0.85) {$C$};
\node at (53:0.91) {$M$};

\node at (60:0.35) {\small $a_1$};
\node at (-50:0.3) {\small $a_2$};
\node at (6:0.75) {\small $b_1$};
\node at (-20:1.12) {\small $b_2$};
\node at (-63:0.57) {\small $c_1$};
\node at (-50:0.98) {\small $c_2$};

\end{tikzpicture}
\caption{Pentagonal subdivision of a regular triangle.}
\label{div_triangle}
\end{figure}

In the second of Figure \ref{div_triangle}, $A,B,C$ are the center of $F$, a vertex of $F$, and the middle point of an edge of $F$. Moreover, let $B'$ and $M$ be the rotations of $B$ and $C$ around the center $A$ by $\frac{2}{3}\pi$. Then $M$ is the middle point of $BB'$.

\begin{figure}[htp]
\centering
\begin{tikzpicture}[>=latex,scale=5]

\node at (-0.2,-1.2) 
 {\includegraphics[scale=0.20]{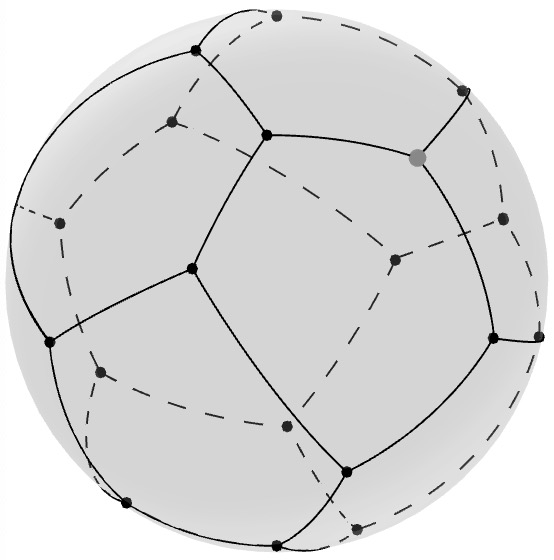}};
\node at (0.06,-1.03){\small $V$};

\node at (0.7,-1.2) 
 {\includegraphics[scale=0.18]{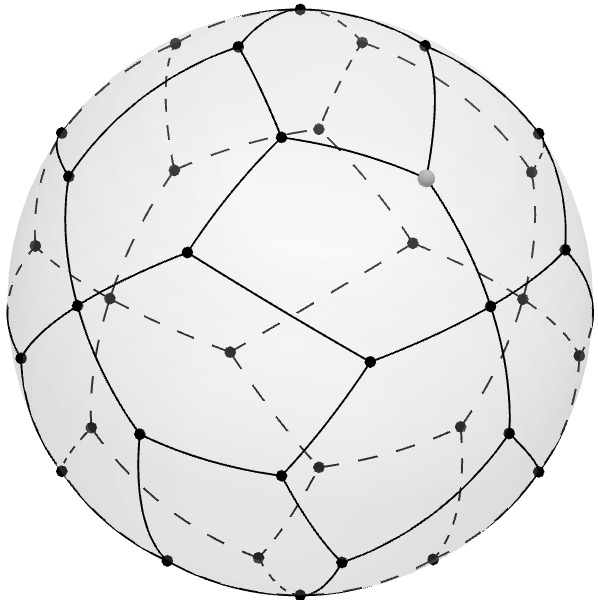}};
\node at (0.91,-1.03){\small $V$};
 
\node at (1.6,-1.2) 
 {\includegraphics[scale=0.17]{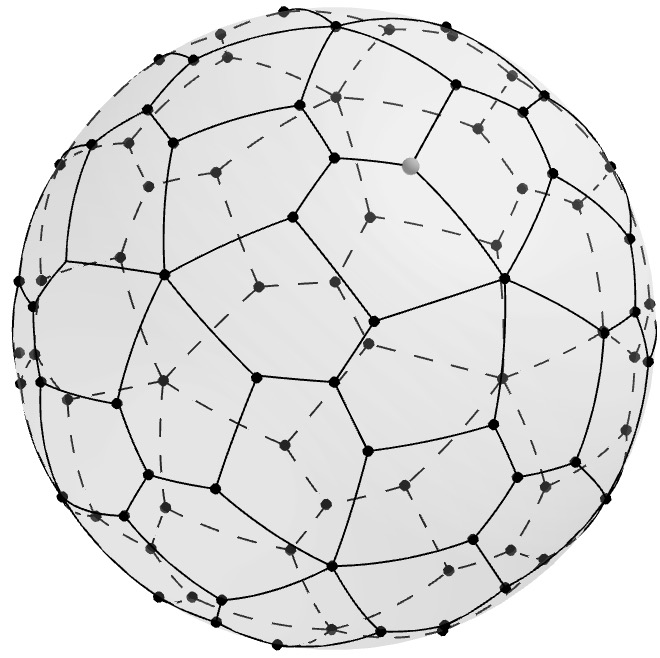}};
\node at (1.75,-1){\small $V$};

\end{tikzpicture}
\caption{Pentagonal subdivision tilings with $12$, $24$, $60$ tiles.}
\label{pst}
\end{figure}

The construction of the pentagonal subdivision starts with an {\em anchor point} $V$ indicated by ${\color{gray} \bullet}$. We connect $V$ to $A$ by the arc of length $<\pi$ (the other arc has length $>\pi$) to get $a_1$, and then rotate $a_1$ around $A$ by $-\frac{2}{3}\pi$ to get $a_2$. On the other side, we connect $V$ to $B$ by the arc of length $<\pi$ to get $b_1$, and then rotate $b_1$ around $B$ by $\frac{2}{n}\pi$ to get $b_2$. We further connect the other ends of $a_2$ and $b_2$ by arcs of length $<\pi$  to $C$. It turns out the two arcs form one arc $c$, and $C$ is the middle point of $c$.

A tiling of the sphere by congruent polygons requires all the tiles to be simple, in the sense that the boundary of each tile is a simple closed curve \cite[Lemma 1]{gsy}. This is the reason for the $<\pi$ requirement in the construction above. Therefore the moduli space is the set of all $V$, such that the corresponding pentagon is simple.

\begin{theorem*}
The moduli space of pentagonal subdivision tilings is given by the locations of the anchor point in an open region of the sphere bounded by two arcs and three curves in Figure \ref{moduliM}. 
\end{theorem*}

Each of the three curves can be interpreted in two ways. The first is the intersection of a hyperbolic cylinder with the sphere. The second is the intersection of a quadratic cone (i.e., quadratic homogeneous equation) with the sphere. See the remark after \eqref{eqB}. Their explicit formulae are obtained by simple transformations of the general form given by \eqref{eqC} through \eqref{eqG}. 

In Section \ref{region}, we divide the sphere into a number of regions, and discuss which regions are completely inside or outside the moduli space. The remaining regions only have some parts contained in the moduli space. In Section \ref{boundary}, we derive the formulae of the boundary of the moduli space parts in these remaining regions, in preferred stereographic projections (called $A$- and $B$-projections). In Section \ref{picture}, we draw perspective pictures of the moduli space in the preferred stereographic $M$-projection, by using M\"obius transforms to translate data in Section \ref{boundary}. In the final Section \ref{feature}, we discuss various features of the moduli space, including the upper bounds of the edges, the reductions by the equalities of edge lengths, and the size of the moduli space. We explicitly express the areas of the three moduli spaces in terms of elliptic integrals. They are roughly 21.5\%, 11.5\%, 4.9\% of the total area of the sphere.

All calculations are done by symbolic calculation software. All pictures are either real 3D pictures or of the perspective
view (except the first of Figure \ref{div_triangle}, and Figure \ref{same_curve}), with curves drawn in parameterized form. All non-decimal values are precise. The decimal values such as $R=0.53308$ mean $0.53308\le R<0.53309$. The only exception is the decimal values in Table \ref{supremum}, where $\sup a=0.6259072384\pi$ mean $0.6259072383\pi<\sup a<0.6259072384\pi$. The obvious reason for the exception is that the values in the table must be upper bounds.

\section{Moduli Space  Region}
\label{region}

We first determine the rough location of the anchor point $V$. The triangle $\triangle ABC$ in Figure \ref{div_triangle} is one sixth of one regular triangular face of the Platonic solid. The angles at $A,B,C$ are respectively $\frac{1}{3}\pi,\frac{1}{n}\pi,\frac{1}{2}\pi$. We rotate the great circle $\bigcirc AB$ around $A$ by multiples of $\frac{1}{3}\pi$ to divide the sphere into six regions. We also rotate the great circle $\bigcirc AB$ around $B$ by multiples of $\frac{1}{n}\pi$ to divide the sphere into $2n$ regions. The intersections of the two divisions divide the sphere into $6n$ regions. 

We use stereographic projections (see Section \ref{sprojection}) to describe what happens on the sphere. A stereographic projection is determined by its origin $X$ and the direction of the real axis, which we call {\em $X$-projection}. We denote the antipode of $P$ by $P^*$. In our pictures, all circles (and straight lines) are {\em arcs}, which are parts of great circles on the sphere. In a stereographic projection, we may visually determine arcs by the fact that arcs passing through $P$ are exactly circles passing through both $P$ and $P^*$.

We draw the division of the sphere in both $A$-projection and $B$-projection. The former is convenient for visualising the rotation from $a_1$ to $a_2$, and the later is convenient for visualising the rotation from $b_1$ to $b_2$. Figure \ref{projection-tetra} gives the perspective pictures of the two stereographic projections for the regular tetrahedron. We label the regions by $\Omega_1,\Omega_2,\dots,\Omega_{18}$. We also indicate $A,B,C,M$ and their antipodes.

\begin{figure}[htp]
\centering
\begin{tikzpicture}[>=latex,scale=1.2]

%% A-projection

\foreach \a in {1,-1}
\draw[yscale=\a]
	(120:1) circle ({sqrt(3)})
	(120:3.3) -- (-60:3);
	
\draw
	(-2.6,0) -- (1.6,0);
	
\fill 
	(1,0) circle (0.05);
\filldraw[fill=white]
	(0,0) circle (0.05);
	
\node at (-0.25,0.15) {$A$};
\node at (1.2,0.15) {$B$};
\node at (-2.3,0.15) {$B^*$};
\node at (-47:0.9) {$C$};
\node at (115:2.85) {$C^*$};
\node at (50:0.95) {$M$};
\node at (-116:2.9) {$M^*$};

\node at (25:0.5) {$\Omega_1$};
\node at (90:0.5) {$\Omega_3$};
\node at (160:1) {$\Omega_5$};

\node at (40:1.3) {$\Omega_7$};
\node at (90:1.5) {$\Omega_9$};
\node at (135:2) {$\Omega_{11}$};

\node at (45:2.1) {$\Omega_{13}$};
\node at (90:2.8) {$\Omega_{15}$};
\node at (140:3) {$\Omega_{17}$};

\node at (-25:0.5) {$\Omega_2$};
\node at (-90:0.5) {$\Omega_4$};
\node at (-160:1) {$\Omega_6$};

\node at (-40:1.3) {$\Omega_8$};
\node at (-90:1.5) {$\Omega_{10}$};
\node at (-135:2) {$\Omega_{12}$};

\node at (-45:2.1) {$\Omega_{14}$};
\node at (-90:2.8) {$\Omega_{16}$};
\node at (-140:3) {$\Omega_{18}$};

%% B-projection

\begin{scope}[xshift=4cm]

\foreach \a in {1,-1}
\draw[yscale=\a]
	(60:1) circle ({sqrt(3)})
	(60:3.3) -- (-120:3);
	
\draw
	(2.6,0) -- (-1.6,0);
	
\fill 
	(0,0) circle (0.05);
\filldraw[fill=white]
	(-1,0) circle (0.05);
	
\node at (0.25,0.15) {$B$};
\node at (-1.25,0.15) {$A$};
\node at (2.25,0.15) {$A^*$};
\node at (230:0.9) {$C$};
\node at (63:2.95) {$C^*$};
\node at (132:0.95) {$M$};
\node at (-65:2.95) {$M^*$};

\node at (155:0.5) {$\Omega_1$};
\node at (140:1.3) {$\Omega_3$};
\node at (135:2.1) {$\Omega_5$};

\node at (90:0.5) {$\Omega_7$};
\node at (90:1.5) {$\Omega_9$};
\node at (90:2.8) {$\Omega_{11}$};

\node at (20:1) {$\Omega_{13}$};
\node at (35:2) {$\Omega_{15}$};
\node at (40:3) {$\Omega_{17}$};

\node at (-155:0.5) {$\Omega_2$};
\node at (-140:1.3) {$\Omega_4$};
\node at (-135:2.1) {$\Omega_6$};

\node at (-90:0.5) {$\Omega_8$};
\node at (-90:1.5) {$\Omega_{10}$};
\node at (-90:2.8) {$\Omega_{12}$};

\node at (-20:1) {$\Omega_{14}$};
\node at (-35:2) {$\Omega_{16}$};
\node at (-40:3) {$\Omega_{18}$};

\end{scope}

\end{tikzpicture}
\caption{$A$- and $B$-projections of $18$ regions for the tetrahedron.}
\label{projection-tetra}
\end{figure}

Figure \ref{core_region} shows the $A$-projection of the pentagon when $V$ lies in $\Omega_1$, $\Omega_2$, $\Omega_3$, $\Omega_7$. The first is the case $V\in \Omega_1$. We have $a_1,b_1\in \Omega_1$. This implies $a_2,c_1\in \Omega_4$ and $b_2,c_2\in\Omega_8$. The locations of these arcs imply that the pentagon is simple. Similarly, for $V\in \Omega_2$, we have $a_1,b_1\in \Omega_2$, $a_2,c_1\in \Omega_4\cup\Omega_6$, $b_2,c_2\in \Omega_8\cup\Omega_{14}$. For $V\in \Omega_3$, we have $a_1,b_1\in \Omega_1\cup\Omega_3$, $a_2,c_1\in \Omega_2$, $b_2,c_2\in \Omega_8\cup\Omega_{10}$. For $V\in \Omega_7$, we have $a_1,b_1\in \Omega_1\cup\Omega_7$, $a_2,c_1\in \Omega_4\cup\Omega_{10}$, $b_2,c_2\in \Omega_2$. All imply that the pentagon is simple. Two particular such pentagonal subdivision tilings are shown in the first two 3D pictures of Figure \ref{pst2}. 

We also consider the extreme case that $V$ lies on the boundaries of these regions. The pentagon is not simple exactly when $V$ is any of the vertices of $\Omega_1,\Omega_2,\Omega_3,\Omega_7$, or $V$ is in the edges $\Omega_3\cap \Omega_9$ ($b_2,c_2$ overlap) or $\Omega_7\cap \Omega_9$ ($a_2,c_1$ overlap).

\begin{figure}[htp]
\centering
\begin{tikzpicture}[>=latex,scale=1]

%% backbround

\foreach \x in {0,1,2,3}
\foreach \a in {0,1,2}
\draw[xshift=3.3*\x cm, rotate=120*\a, gray]
	(1,0) arc (30:120:{sqrt(3)})
	(1,0) arc (30:0:{sqrt(3)})
	(120:1.6) -- (-60:1.3);

\node at (0,0.25) {\small $A$};
\node at (1.2,0) {\small $B$};
\node at (-70:0.9) {\small $C$};
\node at (38:0.83) {\small $V$};

%%% 1

\foreach \a in {0,1,2}
{
\begin{scope}[rotate=120*\a]
	
\coordinate (M) at (60:{sqrt(3)-1});
\coordinate (X) at (0:1);

\pgfmathsetmacro{\th}{18} 

\pgfmathsetmacro{\ph}{-atan((1/3)*(cot(\th+30)))}
\coordinate (A1) at ({\ph}:{sqrt(2)});
\coordinate (B1) at ({\ph+180}:{sqrt(2)});

\pgfmathsetmacro{\ps}{-atan((1/3)*(cot(\th+150)))}
\coordinate (A2) at ({\ps}:{sqrt(2)});
\coordinate (B2) at ({\ps+180}:{sqrt(2)});

\coordinate (V) at (36.8:0.595);
\coordinate (W) at (-42:0.95);
\coordinate (W2) at (78:0.95);

\arcThroughThreePoints[thick]{X}{A1}{V};
\arcThroughThreePoints[thick]{W}{A2}{X};
\arcThroughThreePoints[dashed]{V}{M}{W2};

\draw (0,0) -- (V);

\end{scope}
}

\filldraw[fill=gray!50]
	(36.8:0.595) circle (0.05);

%%% 2

\begin{scope}[xshift=3.3cm]

\foreach \a in {0,1,2}
{
\begin{scope}[rotate=120*\a]

\coordinate (M) at (60:{sqrt(3)-1});
\coordinate (X) at (0:1);
 
\pgfmathsetmacro{\th}{100} 

\pgfmathsetmacro{\ph}{-atan((1/3)*(cot(\th+30)))}
\coordinate (A1) at ({\ph}:{sqrt(2)});
\coordinate (B1) at ({\ph+180}:{sqrt(2)});

\pgfmathsetmacro{\ps}{-atan((1/3)*(cot(\th+150)))}
\coordinate (A2) at ({\ps}:{sqrt(2)});
\coordinate (B2) at ({\ps+180}:{sqrt(2)});

\coordinate (V) at (-26.7:0.64);
\coordinate (W1) at (-10.4:1.7);
\coordinate (W2) at (109.6:1.7);

\arcThroughThreePoints[thick]{W1}{A2}{X};
\arcThroughThreePoints[thick]{V}{A1}{X};
\arcThroughThreePoints[dashed]{V}{M}{W2};

\draw (0,0) -- (V);

\end{scope}
}

\filldraw[fill=gray!50]
	(-26.7:0.64) circle (0.05);

\end{scope}

%%% 3

\begin{scope}[xshift=6.6cm,yscale=-1]
   
\foreach \a in {0,1,2}
{
\begin{scope}[rotate=120*\a]

\coordinate (M) at (60:{sqrt(3)-1});
\coordinate (X) at (0:1);

\pgfmathsetmacro{\th}{-12}

\pgfmathsetmacro{\ph}{-atan((1/3)*(cot(\th+30)))}
\coordinate (A1) at ({\ph}:{sqrt(2)});
\coordinate (B1) at ({\ph+180}:{sqrt(2)});

\pgfmathsetmacro{\ps}{-atan((1/3)*(cot(\th+150)))}
\coordinate (A2) at ({\ps}:{sqrt(2)});
\coordinate (B2) at ({\ps+180}:{sqrt(2)});

\coordinate (V) at (36.8:0.595);
\coordinate (V2) at (-83.2:0.595);
\coordinate (W2) at (78:0.95);

\arcThroughThreePoints[thick]{V2}{A2}{X};
\arcThroughThreePoints[thick]{X}{A1}{W2};
\arcThroughThreePoints[dashed]{V}{M}{W2};

\draw (0,0) -- (V);

\end{scope}
}

\filldraw[fill=gray!50]
	(-83.2:0.595) circle (0.05);

\end{scope}

%%% 4

\begin{scope}[xshift=9.9 cm,yscale=-1]

\foreach \a in {0,1,2}
{
\begin{scope}[rotate=120*\a]

\coordinate (M) at (60:{sqrt(3)-1});
\coordinate (X) at (0:1);

\pgfmathsetmacro{\th}{18} 

\pgfmathsetmacro{\ph}{-atan((1/3)*(cot(\th+30)))}
\coordinate (A1) at ({\ph}:{sqrt(2)});
\coordinate (B1) at ({\ph+180}:{sqrt(2)});

\pgfmathsetmacro{\ps}{-atan((1/3)*(cot(\th+150)))}
\coordinate (A2) at ({\ps}:{sqrt(2)});
\coordinate (B2) at ({\ps+180}:{sqrt(2)});

\coordinate (V) at (36.8:0.595);
\coordinate (W) at (-42:0.95);
\coordinate (W2) at (78:0.95);

\arcThroughThreePoints[thick]{X}{A1}{V};
\arcThroughThreePoints[thick]{W}{A2}{X};
\arcThroughThreePoints[dashed]{V}{M}{W2};

\draw (0,0) -- (78:0.95);

\end{scope}
}

\filldraw[fill=gray!50]
	(-42:0.95) circle (0.05);

\end{scope}

%% backbround

\foreach \x in {0,1,2,3}
{
\begin{scope}[xshift=3.3*\x cm]

\fill (1,0) circle (0.05);
\filldraw[fill=white] (0,0) circle (0.05);	

\end{scope}	
}

\end{tikzpicture}
\caption{The pentagon is simple when $V$ is in $\Omega_1,\Omega_2,\Omega_3,\Omega_7$.}
\label{core_region}
\end{figure}

\begin{figure}[htp]
\centering
\begin{tikzpicture}[>=latex,scale=5]

\node at (-0.2,-1.2) 
 {\includegraphics[scale=0.25]{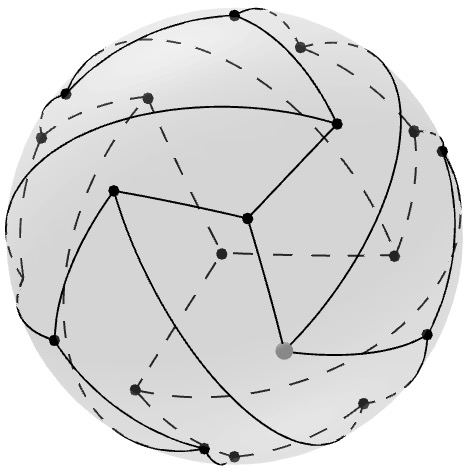}};
\node at (-0.162,-1.39){\small$V$};

\node at (0.7,-1.2) 
 {\includegraphics[scale=0.19]{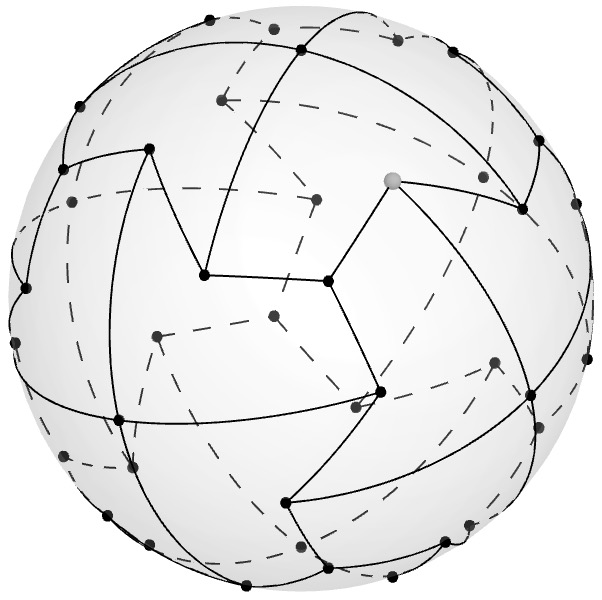}};
  \node at (0.80,-1.02){\small$V$};

\node at (1.6,-1.2) 
 {\includegraphics[scale=0.11]{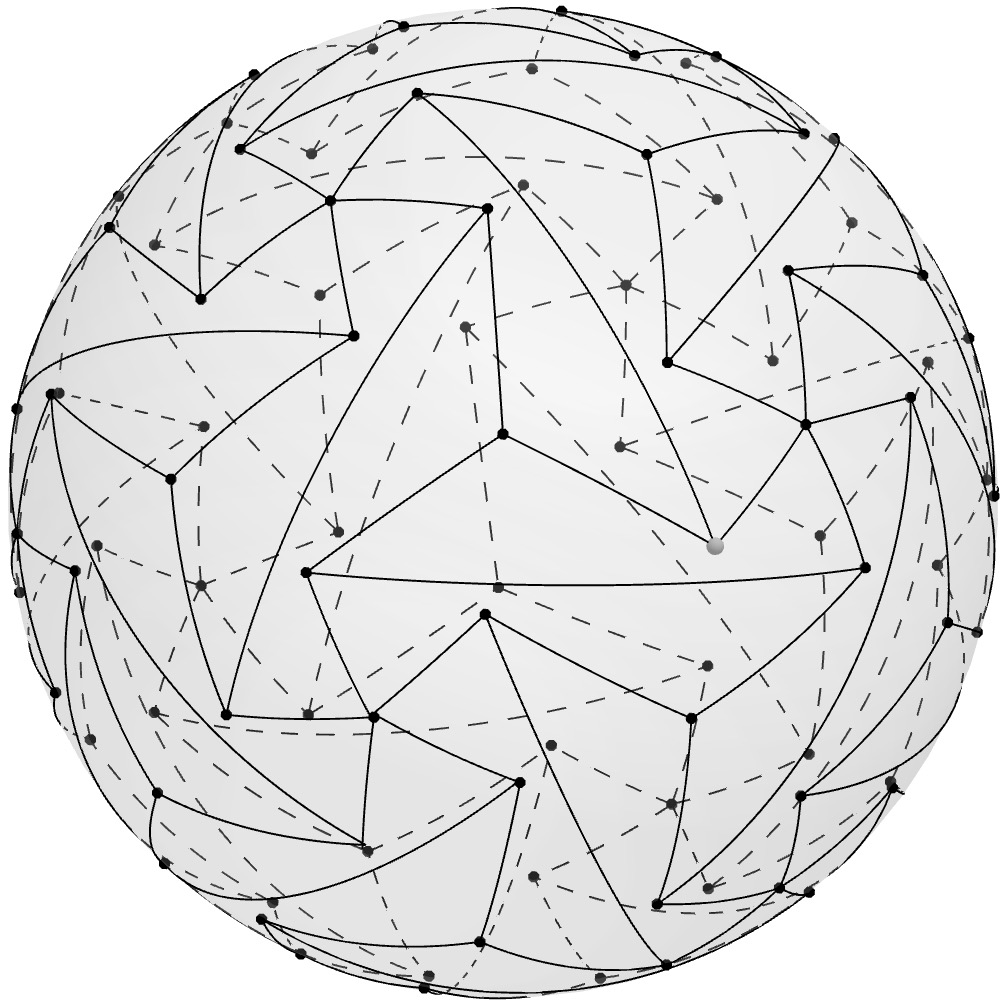}};
 \node at (1.78,-1.18){\small$V$};

\end{tikzpicture}
\caption{Pentagonal subdivision tilings when $V$ is in $\Omega_2,\Omega_3$, or part of $\Omega_8$.}
\label{pst2}
\end{figure}

Next we show that the pentagon is not simple when $V$ is in some other regions. If $V\in \Omega_9$, then $a_1\in A-3-9$, which means that $a_1$ starts from $A$ and then successively passes through $\Omega_3,\Omega_9$. Using the $A$-projection in Figure \ref{projection-tetra}, we may deduce $a_2\in A-2-8$ or $a_2\in A-2-8-14$. We abbreviate this by writing $a_2\in A-2-8-\cdots$. Similarly, we have $b_1\in B-7-9$. Using the $B$-projection, we get $b_2\in B-2-4-\cdots$. Since a curve in the direction $A-2-8$ and a curve in the direction $B-2-4$ must intersect (the intersection point is inside $\Omega_2$), we conclude that $a_2,b_2$ intersect, and the pentagon is not simple.

Similar argument shows that the pentagon is not simple in the following cases. We indicate the two intersecting edges.
\begin{itemize}
\item $V\in \Omega_9$: $a_2\in A-2-8-\cdots$, $b_2\in B-2-4-\cdots$.
\item $V\in \Omega_{11}$: $a_2\in A-1-7-13$, $b_1\in B-7-9-11$.
\item $V\in \Omega_{15}$: $a_1\in A-3-9-15$, $b_2\in B-1-3-5$.
\item $V\in \Omega_{17}$: $a_2\in A-1-7-13$, $b_2\in B-1-3-5$.
\item $V\in \Omega_6$: $b_1\in B-2-4-6$, $c_1\in C-2-1-\cdots$. 
\item $V\in \Omega_{10}$: $a_1\in A-4-10$, $c_1\in C-4-6-5-11$.
\item $V\in \Omega_{12}$: $b_1\in B-8-10-12$, $c_2\in C-10-16-17$.
\item $V\in \Omega_{14}$: $a_1\in A-2-8-14$, $c_2\in C-2-1-(3/7)-9$. Here $(3/7)$ means the arc passes through $\Omega_3$ or $\Omega_7$.
\item $V\in \Omega_{16}$: $a_1\in A-4-10-16$, $c_1\in C-10-12-18-17$.
\item $V\in \Omega_{18}$: $a_2\in A-3-9-15$, $b_2\in B-7-9-11$.
\end{itemize}

For the pentagonal subdivision of the tetrahedron, it remains to study the regions $\Omega_4,\Omega_5,\Omega_8,\Omega_{13}$. Figure \ref{boundary_region} shows that it is possible for the pentagon to be simple if $V$ is in certain parts of these regions. The third of Figure \ref{pst2} shows the 3D picture of one such pentagonal subdivision tiling. We calculate these regions in Section \ref{boundary}.

\begin{figure}[htp]
\centering
\begin{tikzpicture}[>=latex,scale=1]

%% backbround

\foreach \x/\y in {0/1.6, 3.6/1.6, 1.8/-1.2}
\foreach \a in {0,1,2}
\draw[shift={(\x cm,\y cm)}, rotate=120*\a, gray]
	(1,0) arc (30:120:{sqrt(3)})
	(1,0) arc (30:0:{sqrt(3)})
	(120:1.6) -- (-60:1.3);
	
%%% 1

\begin{scope}[shift={(0cm,1.6cm)}]

\foreach \a in {0,1,2}
{
\begin{scope}[rotate=120*\a]

\coordinate (B) at (1,0);
\coordinate (M1) at (-60:{sqrt(3)-1});
\coordinate (M2) at (120:{sqrt(3)+1});
\coordinate (V) at (130-120:0.71);
\coordinate (W) at (-109:1.33);

\draw (0,0) -- (130:0.71);

\draw[thick]
	(B) arc (50:89:2.333)
	(B) arc (-10:-87.5:1.523);
	
\arcThroughThreePoints[dashed]{W}{M1}{V};

\end{scope}
}

\fill (1,0) circle (0.05);
\filldraw[fill=white] (0,0) circle (0.05);	

\filldraw[fill=gray!50]
	 (129:0.7) circle (0.05);

\end{scope}

%%% 2

\begin{scope}[shift={(3.6cm,1.6cm)}]

\foreach \a in {0,1,2}
{
\begin{scope}[rotate=120*\a]

\coordinate (B1) at (1,0);
\coordinate (B2) at (-2,0);
\coordinate (M1) at (-60:{sqrt(3)-1});
\coordinate (M2) at (120:{sqrt(3)+1});

\coordinate (P) at (6.2:0.579);

\coordinate (V) at (20:1.376);
\coordinate (W) at (20-120:1.376);

\arcThroughThreePoints[thick]{B1}{B2}{V};
\arcThroughThreePoints[thick]{B1}{B2}{P};
\arcThroughThreePoints[dashed]{W}{M1}{P};
	
\draw
	(0,0) -- (V);

\end{scope}
}

\fill (1,0) circle (0.05);
\filldraw[fill=white] (0,0) circle (0.05);	

\filldraw[fill=gray!50]
	 (20:1.376) circle (0.05);

\end{scope}

%%% 3

\begin{scope}[shift={(1.8cm,-1.2cm)}]

\foreach \a in {0,1,2}
{
\begin{scope}[rotate=120*\a]

\coordinate (B1) at (1,0);
\coordinate (B2) at (-2,0);
\coordinate (M1) at (-60:{sqrt(3)-1});
\coordinate (M2) at (120:{sqrt(3)+1});

\coordinate (P) at (8:1.493);

\coordinate (V) at (-20:0.904);
\coordinate (W) at (-20-120:0.904);

\arcThroughThreePoints[thick]{V}{B2}{B1};
\arcThroughThreePoints[thick]{B1}{B2}{P};
\arcThroughThreePoints[dashed]{W}{M1}{P};	
	
\draw
	(0,0) -- (V);

\end{scope}
}

\fill (1,0) circle (0.05);
\filldraw[fill=white] (0,0) circle (0.05);	

\filldraw[fill=gray!50]
	 (-20:0.904) circle (0.05);

\end{scope}

%%% 4

\begin{scope}[xshift=7.8cm]

\foreach \a in {0,1,2} 
\draw[rotate=120*\a, gray]
	(1,0) arc (30:390:{sqrt(3)})
	(-2.3,0) -- (3,0);

\foreach \a in {0,1,2}
{
\begin{scope}[rotate=120*\a]

\coordinate (B1) at (1,0);
\coordinate (B2) at (-2,0);
\coordinate (M1) at (-60:{sqrt(3)-1});
\coordinate (M2) at (120:{sqrt(3)+1});

\coordinate (P) at (-33.37:2.75);

\coordinate (V) at (-70:0.358);
\coordinate (W) at (-70-120:0.358);

\arcThroughThreePoints[thick]{P}{B2}{B1};
\arcThroughThreePoints[thick]{V}{B2}{B1};
\arcThroughThreePoints[dashed]{W}{M1}{P};

\draw
	(0,0) -- (V);

\end{scope}
}

\filldraw[fill=gray!50]
	 (-66:0.35) circle (0.05);

\fill (1,0) circle (0.05);
\filldraw[fill=white] (0,0) circle (0.05);	

\end{scope}

\end{tikzpicture}
\caption{The pentagon is simple when $V$ is in certain parts of $\Omega_5,\Omega_{13},\Omega_8,\Omega_4$.}
\label{boundary_region}
\end{figure}

We may carry out similar region by region argument for the pentagonal subdivision of the octahedron. Figure \ref{projection-octa} is the two stereographic projections for the regular octahedron. We simplify the presentation by labelling the regions by $i$ instead of $\Omega_i$. If $V$ is in $\Omega_1,\Omega_2,\Omega_3,\Omega_7$, then the same argument shows that the pentagon is simple. Moreover, the pentagon is not simple if $V$ is a vertex of these regions, or is in $\Omega_3\cap \Omega_9$ or $\Omega_7\cap \Omega_9$.

\begin{figure}[htp]
\centering
\begin{tikzpicture}[>=latex,scale=1]

% A-projection 

\foreach \a in {1,2}
\draw[rotate=120*\a]
	(-3.7,0) -- (4.7,0);

\foreach \a in {0,2}
\draw[rotate=120*\a]
	(-1,{-sqrt(3)}) circle ({sqrt(6)});

\draw
	(-1,0) circle ({sqrt(3)})
	(-3.7,0) -- (2.5,0);	

\fill	
	({sqrt(3)-1},0) circle (0.06);
\filldraw[fill=white]
	(0,0) circle (0.06);
	
\node at (150:0.35) {\small $A$};
\node at (175.5:2.2) {$B^*$};
\node at (9:1.1) {\small $B$};
\node at (-75:0.6) {\scriptsize $C$};
\node at (116:4.2) {$C^*$};
\node at (73:0.63) {\scriptsize $M$};
\node at (-116:4.15) {$M^*$};

\node at (23:0.38) {\scriptsize $1$};
\node at (90:0.32) {\scriptsize $3$};
\node at (160:1) {$5$};

\node at (40:0.67) {\small $7$};
\node at (95:1) {$9$};
\node at (145:1.8) {$11$};

\node at (45:1.4) {$13$};
\node at (95:2.7) {$15$};
\node at (140:3.2) {$17$};

\node at (30:2) {$19$};
\node at (70:3.6) {$21$};
\node at (170:3.6) {$23$};

\node at (-23:0.38) {\scriptsize $2$};
\node at (-90:0.32) {\scriptsize $4$};
\node at (-160:1) {$6$};

\node at (-40:0.67) {\small $8$};
\node at (-95:1) {$10$};
\node at (-145:1.8) {$12$};

\node at (-45:1.4) {$14$};
\node at (-95:2.7) {$16$};
\node at (-140:3.2) {$18$};

\node at (-30:2) {$20$};
\node at (-70:3.6) {$22$};
\node at (-170:3.6) {$24$};

% B-projection 

\begin{scope}[xshift=5cm,scale=1.3]

\foreach \b in {1,-1}
\draw[yscale=\b]
	(1,-1) circle (2);
	
\draw
	(0,-3.2) -- (0,3.2)
	(-1.5,0) -- (3.2,0)
	(-1.5,-1.5) -- (3,3)
	(-1.5,1.5) -- (3,-3);

\fill
	(0,0) circle (0.046);
	
\filldraw[fill=white]
	({-sqrt(3)+1},0) circle (0.046);
	
\node at (170:1) { $A$};
\node at (3,0.15) {$A^*$};
\node at (25:0.4) {$B$};
\node at (235:0.8) {$C$};
\node at (2.4,2.15) {$C^*$};
\node at (120:0.8) {$M$};
\node at (2.45,-2.1) {$M^*$};

\node at (157:0.4) {\small $1$};
\node at (150:0.8) {$3$};
\node at (165:1.3) {$5$};

\node at (115:0.4) {\small $7$};
\node at (110:1.2) {$9$};
\node at (110:2.5) {$11$};

\node at (65:0.65) {$13$};
\node at (70:2) {$15$};
\node at (50:3.7) {$17$};

\node at (20:1.2) {$19$};
\node at (30:2.7) {$21$};
\node at (40:3.7) {$23$};

\node at (-157:0.4) {\small $2$};
\node at (-150:0.8) {$4$};
\node at (-165:1.3) {$6$};

\node at (-115:0.4) {\small $8$};
\node at (-110:1.2) {$10$};
\node at (-110:2.5) {$12$};

\node at (-65:0.65) {$14$};
\node at (-70:2) {$16$};
\node at (-50:3.7) {$18$};

\node at (-20:1.2) {$20$};
\node at (-30:2.7) {$22$};
\node at (-40:3.7) {$24$};

\end{scope}

\end{tikzpicture}
\caption{$A$- and $B$-projections of $24$ regions for the octahedron.}
\label{projection-octa}
\end{figure}

If $V$ is in the following regions, then the pentagon is not simple. Again, we indicate the two intersecting edges.
\begin{itemize}
\item $V\in\Omega_9$: $a_2\in A-2-8-\cdots$, $b_2\in B-2-4-\cdots$.
\item $V\in\Omega_{11}$: $a_2\in A-1-7-13-19$, $b_1\in B-7-9-11$.
\item $V\in\Omega_{15}$: $a_1\in A-3-9-15$, $b_2\in B-1-3-5$.
\item $V\in\Omega_{17}$: $a_2\in A-1-7-13-19$, $b_1\in B-13-15-17$.
\item $V\in\Omega_{21}$: $a_1\in A-3-9-15-21$, $b_2\in B-7-9-11$.
\item $V\in\Omega_{23}$: $a_2\in A-1-7-13-19$, $b_2\in B-7-9-11$.
\item $V\in\Omega_6$: $b_1\in B-2-4-6$, $c_1\in C-2-1-\cdots$.
\item $V\in\Omega_{10}\cup\Omega_{16}$: $a_1\in A-4-10-\cdots$, $c_1\in C-4-6-\cdots$.
\item $V\in\Omega_{20}$: $a_1\in A-2-8-14-20$, $c_2\in B-2-1-\cdots$.
\item $V\in\Omega_{22}$: $a_1\in A-4-10-16-22$, $b_2\in B-8-10-12$.
\item $V\in\Omega_{24}$: $a_2\in A-3-9-15-21$, $b_2\in B-13-15-17$.
\end{itemize}

For some other regions, we need to use the following geometric lemma to show that the pentagon is not simple.

\begin{lemma}\label{geometry1}
Suppose $P$ and $P^*$ are a pair of antipodal points, and $P^*X>PX$ and $VX=WX$ in Figure \ref{geometry}.
\begin{enumerate}
\item If $\rho\le\frac{1}{2}\pi$ and $\theta^*\le\theta$, then $WX$ and $PV$ intersect.
\item If $\rho\ge\frac{1}{2}\pi$ and $\theta^*\ge\theta$, then $WX$ and $PV$ intersect.
\end{enumerate}
\end{lemma}

\begin{figure}[htp]
\centering
\begin{tikzpicture}[>=latex,scale=1]
	
%% 1

\draw
	(-2,0) -- (2,0)
	(0.3,0) -- ++(150:1.06)
	(0.3,0) -- ++(50:1.06)
	(1.5,0) arc (45:135:2.1213);
	
\draw[dashed]
	(1.5,0) arc (0:180:1.5);

\filldraw[fill=white]
	(0.3,0) circle (0.05);
	
\node at (0.3,-0.25) {$X$};
\node at (1.5,-0.2) {$P$};
\node at (-1.5,-0.2) {$P^*$};
\node at (-0.6,0.7) {$V$};
\node at (0.8,0.9) {$W$};
\node at (1.08,0.15) {\small $\rho$};
\node at (-1.12,0.15) {\small $\rho$};
\node at (0.65,0.2) {\small $\theta$};
\node at (-0.4,0.2) {\small $\theta^*$};

%% 2

\begin{scope}[xshift=5cm]

\draw
	(-2,0) -- (2,0)
	(0.3,0) -- ++(130:2.21)
	(0.3,0) -- ++(30:2.21)
	(1.5,0) arc (-20:200:1.5963);
		
\draw[dashed]
	(1.5,0) arc (0:180:1.5);

\filldraw[fill=white]
	(0.3,0) circle (0.05);
	
\node at (0.3,-0.25) {$X$};
\node at (1.5,-0.2) {$P$};
\node at (-1.5,-0.2) {$P^*$};
\node at (-1.2,1.8) {$V$};
\node at (2.45,1.1) {$W$};
\node at (1.37,0.15) {\small $\rho$};
\node at (-1.37,0.15) {\small $\rho$};
\node at (0.9,0.18) {\small $\theta$};
\node at (-0.2,0.2) {\small $\theta^*$};

\end{scope}

\end{tikzpicture}
\caption{A geometric lemma, presented in $X$-projection.}
\label{geometry}
\end{figure}

Figure \ref{geometry} is the $X$-projection. The lines are arcs, the dashed line is a half circle, and all angles are positive. It is easy to see the validity of the lemma in the Euclidean metric of ${\bb C}={\bb R}^2$ ($PXP^*$, $VX$, $WX$ are literally straight lines). We note that, in the $X$-projection, $P^*X>PX$ and $VX=WX$ in the Euclidean metric are equivalent to the same (in)equalities in the spherical metric. Since the stereographic projection preserves the angles, the lemma is actually valid on the sphere, and therefore is also valid in any stereographic projection as long as the lengths are understood to be spherical.

Applying Lemma \ref{geometry1}, the pentagon is not simple in the following cases.
\begin{itemize}
\item $V\in \Omega_{12}$: $X=B$, $P=C$, $\rho\ge\frac{1}{2}\pi$, $\theta^*\ge\frac{1}{4}\pi\ge\theta$, $WX=b_1,PV=c_2$.
\item $V\in \Omega_{14}$: $X=B$, $P=C$, $\rho\le\frac{1}{2}\pi$, $\theta^*\le\frac{1}{4}\pi\le\theta$, $WX=b_1,PV=c_2$.
\item $V\in \Omega_{18}$: $X=B$, $P=C$, $\rho\ge\frac{1}{2}\pi$, $\theta^*\ge\frac{1}{4}\pi\ge\theta$, $WX=b_1,PV=c_2$.
\item $V\in \Omega_{19}$: $X=B$, $P=C$, $\rho\le\frac{1}{2}\pi$, $\theta^*\le\frac{1}{4}\pi\le\theta$, $WX=b_2,PV=a_1$.
\end{itemize}

Again it remains to study the regions $\Omega_4,\Omega_5,\Omega_8,\Omega_{13}$ for the pentagonal subdivision of the octahedron. 

Figures \ref{projection-icoA} and \ref{projection-icoB}  are the two stereographic projections for the icosahedron. If $V$ is in $\Omega_1,\Omega_2,\Omega_3,\Omega_7$, then the same argument shows that the pentagon is simple. Moreover, the pentagon is not simple if $V$ is a vertex of these regions, or is in $\Omega_3\cap \Omega_9$ or $\Omega_7\cap \Omega_9$.

\begin{figure}[htp]
\centering
\begin{tikzpicture}[>=latex,scale=1.7]

% A-projection 

\coordinate (B) at ({(sqrt(30+6*sqrt(5))-sqrt(5)-3)/4},0);
\coordinate (B1) at ({-4/(sqrt(30+6*sqrt(5))-sqrt(5)-3)},0);
\pgfmathsetmacro{\ra}{2.80252}
\pgfmathsetmacro{\rb}{1.73205}

\coordinate (E) at (0.5,-5);
\coordinate (V) at (-1.5,-1.35);

\arcThroughThreePoints[thick]{E}{B1}{B};
\arcThroughThreePoints[thick]{V}{B1}{B};

\draw
	(-4,0) -- (1.8,0)
	(V) -- (0,0);

\foreach \b in {1,-1}
\draw[yscale=\b]
	(60:3.5) -- (240:6)
	(B) arc (54:180:\ra) % low left
	(B) arc (54:-180:\ra) % low right
	(B) arc (18:180:\rb) % high left
	(B) arc (18:-180:\rb); % high right	

\draw[dashed]
	(E) -- (0.095,-0.16);
	
\fill
	(B) circle (0.03);
	
\filldraw[fill=white]
	(0,0) circle (0.03);

\filldraw[fill=gray!50]
	(V) circle (0.05);
	
\node at (-0.18,0.1) {$A$};
\node at (0.61,0.1) {$B$};
\node at (-3.3,0.1) {$B^*$};
\node at (0.03,-0.3) {\small $C$};
\node at (-2.9,4.7) {$C^*$};
\node at (0.03,0.34) {\small $M$};
\node at (-2.9,-4.75) {$M^*$};
\node at (0.35,-5) {$W$};
\node at (-1.6,-1.5) {$V$};

\node at (-0.7,-0.8) {\small $a_1$};
\node at (-0.05,-1) {\small $b_1$};
\node at (1.8,-2) {\small $b_2$};
\node at (0.4,-2.5) {\small $c_2$};

\node at (0.12,0.07) {\scriptsize 1};
\node at (0,0.13) {\scriptsize 3};
\node at (-1.2,0.3) {5};

\node at (0.19,0.18) {\scriptsize 7};
\node at (-0.1,0.5) {9};
\node at (-1.5,0.8) {11};

\node at (0.32,0.35) {\scriptsize 13};
\node at (-0.2,1.2) {15};
\node at (-2,1.6) {17};

\node at (1,1.2) {19};
\node at (-0.4,3.3) {21};
\node at (-3,2.7) {23};

\node at (1.5,0.6) {25};
\node at (1,4.5) {27};
\node at (-3.5,4.5) {29};

\node at (0.12,-0.07) {\scriptsize 2};
\node at (0,-0.13) {\scriptsize 4};
\node at (-1.2,-0.3) {6};

\node at (0.19,-0.18) {\scriptsize 8};
\node at (-0.1,-0.5) {10};
\node at (-1.5,-0.9) {12};

\node at (0.32,-0.35) {\scriptsize 14};
\node at (-0.2,-1.4) {16};
\node at (-2,-1.6) {18};

\node at (1,-1.2) {20};
\node at (-0.4,-3.3) {22};
\node at (-3,-2.7) {24};

\node at (1.5,-0.6) {26};
\node at (-0.2,-5.2) {28};
\node at (-3.7,-4.6) {30};

\end{tikzpicture}
\caption{$A$-projections of $30$ regions for the icosahedron, and the case $V\in \Omega_{18}$.}
\label{projection-icoA}
\end{figure}

\begin{figure}[htp]
\centering
\begin{tikzpicture}[>=latex,scale=2.5]

% B-projection 

\coordinate (A) at ({-(sqrt(30+6*sqrt(5))-sqrt(5)-3)/4},0);
\pgfmathsetmacro{\ra}{sqrt((sqrt(5)+5)/2)};

\draw
	(-0.8,0) -- (3.4,0);

\foreach \b in {1,-1}
\draw[yscale=\b]
	(36:4) -- (216:1)
	(72:3.2) -- (252:2.3)
	(A) arc (150:510:\ra);

\fill
	(0,0) circle (0.02);
	
\filldraw[fill=white]
	(A) circle (0.02);

\node at (-0.47,0.07) {$A$};
\node at (3.1,0.1) {$A^*$};
\node at (0.2,0.07) {$B$};
\node at (-0.26,-0.28) {\small $C$};
\node at (3,2.07) {$C^*$};
\node at (-0.25,0.28) {\small $M$};
\node at (3,-2.03) {$M^*$};

\node at (-0.2,0.07) {\scriptsize 1};
\node at (-0.34,0.15) {\scriptsize 3};
\node at (-0.7,0.2) {5};

\node at (-0.13,0.17) {\scriptsize 7};
\node at (-0.4,0.7) {9};
\node at (-0.7,1.3) {11};

\node at (0,0.28) {\scriptsize 13};
\node at (0,1.4) {15};
\node at (0,2.8) {17};

\node at (0.4,0.5) {19};
\node at (1.3,1.7) {21};
\node at (2.6,2.8) {23};

\node at (1.5,0.4) {25};
\node at (2.6,1) {27};
\node at (3.2,1.7) {29};

\node at (-0.2,-0.07) {\scriptsize 2};
\node at (-0.34,-0.15) {\scriptsize 4};
\node at (-0.7,-0.2) {6};

\node at (-0.13,-0.17) {\scriptsize 8};
\node at (-0.4,-0.7) {10};
\node at (-0.7,-1.3) {12};

\node at (0,-0.28) {\scriptsize 14};
\node at (0,-1.4) {16};
\node at (0,-2.8) {18};

\node at (0.4,-0.5) {20};
\node at (1.3,-1.7) {22};
\node at (2.6,-2.8) {24};

\node at (1.5,-0.4) {26};
\node at (2.6,-1) {28};
\node at (3.2,-1.7) {30};

\end{tikzpicture}
\caption{$B$-projections of $30$ regions for the icosahedron.}
\label{projection-icoB}
\end{figure}

If $V$ is in the following regions, then the pentagon is not simple.
\begin{itemize}
\item $V\in\Omega_9$: $a_2\in A-2-8-\cdots$, $b_2\in B-2-4-\cdots$.
\item $V\in\Omega_{11}$: $a_2\in A-1-7-13-19-25$, $b_1\in B-7-9-11$.
\item $V\in\Omega_{15}$: $a_1\in A-3-9-15$, $b_2\in B-1-3-5$.
\item $V\in\Omega_{17}$: $a_2\in A-1-7-13-19-25$, $b_1\in B-13-15-17$.
\item $V\in\Omega_{21}$: $a_1\in A-3-9-15-21$, $b_2\in B-7-9-11$.
\item $V\in\Omega_{23}$: $a_2\in A-1-7-13-19-25$, $b_2\in B-7-9-11$.
\item $V\in\Omega_{27}$: $a_1\in A-3-9-15-21-27$, $b_2\in B-13-15-17$.
\item $V\in\Omega_{29}$: $a_2\in A-1-7-13-19-25$, $b_2\in B-13-15-17$.
\item $V\in\Omega_6$: $b_1\in B-2-4-6$, $c_1\in C-2-1-\cdots$.
\item $V\in\Omega_{10}\cup\Omega_{16}\cup\Omega_{22}$: $a_1\in A-4-10-\cdots$, $c_1\in C-4-6-\cdots$. 
\item $V\in\Omega_{26}$: $a_1\in A-2-8-14-20-26$, $c_2\in C-2-1-\cdots$.
\item $V\in\Omega_{30}$: $a_2\in A-3-9-15-21-27$, $b_2\in B-19-21-23$.
\end{itemize}

Applying Lemma \ref{geometry1}, the pentagon is not simple in the following cases.
\begin{itemize}
\item $V\in \Omega_{12}$: $X=B$, $P=C$, $\rho\ge\frac{1}{2}\pi$, $\theta^*\ge\frac{2}{5}\pi\ge\theta$, $WX=b_1,PV=c_2$.
\item $V\in \Omega_{20}$: $X=B$, $P=C$, $\rho\le\frac{1}{2}\pi$, $\theta^*\le\frac{2}{5}\pi\le\theta$, $WX=b_1,PV=c_2$.
\item $V\in \Omega_{25}$: $X=B$, $P=A$, $\rho\le\frac{1}{2}\pi$, $\theta^*\le\frac{2}{5}\pi\le\theta$, $WX=b_2,PV=a_1$.
\end{itemize}

For $V\in \Omega_{18},\Omega_{24},\Omega_{29}$, we observe the following.
\begin{itemize}
\item $V\in \Omega_{18}$: $a_1\in A-6-12-18$, $b_1\in B-14-16-18$. The two ends of $b_2$ are $B$ and $W$, with $W\in\Omega_{26}\cup \Omega_{28}\cup \Omega_{30}$. See Figure \ref{projection-icoA}.
\item $V\in\Omega_{24}$: $a_1\in A-6-12-18-24$, $b_1\in B-20-22-24$. The two ends of $b_2$ are $B$ and $W$, with $W\in\Omega_{29}$.
\item $V\in\Omega_{28}$: $a_1\in A-4-10-16-22-28$, $b_1\in B-26-28$. The two ends of $a_2$ are $A$ and $W$, with $W\in\Omega_{29}$.
\end{itemize}
We note that $C$ is inside the triangle $\triangle_C$ bounded by $a_1,b_1$ and $AB$, and $W$ is outside the triangle. Therefore $CW$ intersects the boundary of $\triangle_C$. Since $CW$ is $c_2$ for $V\in \Omega_{18},\Omega_{24}$, and $CW$ is $c_1$ for $V\in \Omega_{28}$, we know $c$ does not intersect $a_1$ and $b_1$. Therefore $CW$ intersect $AB$. 

On the other hand, for the case $V\in \Omega_{18}$, $C$ and $W$ are in the same side of $AB$ (i.e., both points are in the same semisphere divided by the great circle $\bigcirc AB$). For the case $V\in \Omega_{24},\Omega_{28}$, $B$ and $W$ are in different sides of $AC$ (i.e., the two points are in different semispheres divided by the great circle $\bigcirc AC$). Both properties imply that $CW$ does not intersect $AB$.

It remains to study the regions $\Omega_4,\Omega_5,\Omega_8,\Omega_{13},\Omega_{14},\Omega_{19}$ for the pentagonal subdivision of the icosahedron.

\section{Boundary of Moduli Space}
\label{boundary}

For tetrahedron and octahedron, it remains to consider $V$ in $\Omega_4,\Omega_5,\Omega_8,\Omega_{13}$. For icosahedron, it remains to consider $V$ in $\Omega_4,\Omega_5,\Omega_8,\Omega_{13},\Omega_{14},\Omega_{19}$. The anchor point $V$ must be in certain parts of these regions for the pentagon to be simple. The parts are calculated in suitable stereographic projections.

\subsection{Stereographic Projection}
\label{sprojection}

We parameterise a stereographic projection by complex numbers, such that the equator is the circle $z\bar{z}=1$ of radius $1$ centered at $0$. See Figure \ref{projection}. Under the projection, a point $\xi=(\xi_{12},\xi_3)=(\xi_1,\xi_2,\xi_3)$ on the sphere and its complex paramaterisation $z=re^{i\theta}$ are related by
\[
\xi=\frac{(2re^{i\theta},r^2-1)}{r^2+1}
=\frac{(2z,|z|^2-1)}{|z|^2+1},\quad
z=\frac{\xi_{12}}{1-\xi_3},\quad
|z|^2=\frac{1+\xi_3}{1-\xi_3}.
\]

For the convenience of writing the formulae for the boundaries of the moduli spaces in the later discussion, we introduce $R$ by (see \eqref{eqE} and \eqref{eqD})
\begin{equation}\label{eqR1}
r^{-1}-r=2R,\quad 
r=\sqrt{R^2+1}-R.
\end{equation}
Then we have
\begin{equation}\label{eqR2}
\xi=\frac{(e^{i\theta},-R)}{\sqrt{R^2+1}}=\frac{(\cos\theta,\sin\theta,-R)}{\sqrt{R^2+1}}.
\end{equation}

\begin{figure}[htp]
\centering
\begin{tikzpicture}[>=latex,scale=1.5]

\draw
	(0,0) circle (1)
	(-1.5,0) -- (1.5,0)
	(0,1) -- (-35:1);

\draw[gray]
	(0,-1) -- (0,1);

\draw[dotted]
	(0,0) -- (-35:1); 

\node at (1.1,-0.15) {$1$};
\node at (0.6,0.1) {$z$};
\node at (0.9,-0.65) {$\xi$};
\node at (0,-1.15) {$X$};
\node at (-0.1,-0.15) {$0$};
\node at (1.65,0) {${\bb C}$};
\node at (-0.8,0.8) {${\bb S}^2$};
\node at (0.12,-0.25) {$\delta$};

\end{tikzpicture}
\caption{$X$-projection.}
\label{projection}
\end{figure}

The antipode $z^*$ of $z$ is given by $z^*\bar{z}=-1$. The spherical distance $\delta$ between $X$ and $\xi$ satisfies $\tan\frac{1}{2}\delta=|z|$. The parameterisations of different projections are related by M\"obius transforms.

We use $z_X$ to denote the parameterisation by the $X$-projection. For example, we write a M\"obius transform from $X$-projection to $Y$-projection as
\[
z_Y=\frac{az_X+b}{cz_X+d}.
\]
Moreover, we denote by $Y_X\in {\bb C}$ the $X$-projection of $Y\in {\bb S}^2$. Then we always have $X_X=0$ and $X^*_X=\infty$. The geometry shows that the Euclidean distance $|Y_X|$ (norm of complex number) of $X$ and $Y$ in the $X$-projection is equal to the Euclidean distance $|X_Y|$ of $X$ and $Y$ in the $Y$-projection. We denote $d_{XY}=|Y_X|=|X_Y|$.

Our description of the moduli space is based on the triangle $\triangle ABM$, which is the region $\Omega_1$ in Section \ref{region}. The triangle has angles $\frac{1}{3}\pi,\frac{1}{n}\pi,\frac{1}{2}\pi$ at $A,B,M$. We will use the following projections:
\begin{itemize}
\item $A$-projection: $B$ is in the positive real direction.
\item $B$-projection: $A$ is in the negative real direction.
\item $M$-projection: $A$ is in the negative real direction.
\end{itemize} 
Moreover, the direction $A\to B\to M$ is always counterclockwise. This implies that $M$ is in the upper half of the $A$- and $B$-projections, and $B$ is in the lower half (in fact the negative imaginary direction) of the $M$-projection.

The first of Figure \ref{ABM} is the triangle $\triangle ABM$ in the $A$-projection. We have
\[
B_A=d_{AB},\;
B^*_A=-\bar{B}_A^{-1}=-d_{AB}^{-1},\;
M_A=d_{AM}e^{i\frac{1}{3}\pi},
\]
and $BMB^*$ is a circle centered at $O$. Then the middle point $D$ of $B$ and $B^*$ has $D_A=\frac{1}{2}(d_{AB}-d_{AB}^{-1})$. This gives the Euclidean distances in the $A$-projection (indicated by subscript $A$) 
\[
|AD|_A=\tfrac{1}{2}(d_{AB}^{-1}-d_{AB}),\quad
|BD|_A=\tfrac{1}{2}(d_{AB}^{-1}+d_{AB}).
\]
By $\angle OAD=\frac{1}{3}\pi$ and $\angle OBD=\frac{1}{2}\pi-\angle ABM=(\frac{1}{2}-\frac{1}{n})\pi$, we may calculate the distance $|OD|_A$ in two ways and get an equation
\[
(d_{AB}^{-1}-d_{AB})\tan\tfrac{1}{3}\pi
=(d_{AB}^{-1}+d_{AB})\tan(\tfrac{1}{2}-\tfrac{1}{n})\pi.
\]
We solve the equation and get one positive solution  
\[
d_{AB}=\sqrt{\frac{\sqrt{3}\tan\frac{1}{n}\pi-1}{\sqrt{3}\tan\frac{1}{n}\pi+1}}.
\]
The explicit values are given in Table \ref{distance}. Then we further get
\begin{align*}
d_{AM}
&=|OM|_A-|OA|_A=|OB|_A-|OA|_A \\
&=|BD|_A\sec(\tfrac{1}{2}-\tfrac{1}{n})\pi-|AD|_A\sec\tfrac{1}{3}\pi  \\
&=\tfrac{1}{2}(d_{AB}^{-1}+d_{AB})\sec(\tfrac{1}{2}-\tfrac{1}{n})\pi
-\tfrac{1}{2}(d_{AB}^{-1}-d_{AB})\sec\tfrac{1}{3}\pi. 
\end{align*}

\begin{table}[htp]
\centering
\begin{tabular}{|c||c|c|c|}
\hline 
& $n=3$   
& $n=4$  
& $n=5$  \\
& tetrahedron
& octahedron 
& icosahedron \\
\hline  \hline 
$d_{AB}$
& $\frac{1}{\sqrt{2}}$
& $\frac{\sqrt{3}-1}{\sqrt{2}}$
& $\frac{1}{4}(\sqrt{30+6\sqrt{5}}-\sqrt{5}-3)$ \\
\hline
$d_{AM}$
& $\frac{\sqrt{3}-1}{\sqrt{2}}$
& $\sqrt{3}-\sqrt{2}$
& $\frac{1}{2}(\sqrt{15}-\sqrt{5}+\sqrt{3}-3)$ \\
\hline
$d_{BM}$
& $\frac{\sqrt{3}-1}{\sqrt{2}}$
& $\sqrt{2}-1$
& $\frac{1}{2}(\sqrt{10+2\sqrt{5}}-\sqrt{5}-1)$ \\
\hline 
\end{tabular}
\caption{Euclidean distances.}
\label{distance}
\end{table}

\begin{figure}[htp]
\centering
\begin{tikzpicture}[>=latex,scale=1]

\foreach \a in {1,-1}
\draw[xshift=-3*\a cm, xscale=\a]
	(-2.5,{sqrt(3)}) -- (2.5,{sqrt(3)})
	(0,0) -- (90:2.7)
	(0,0) -- (60:2.7)
	(0,0) -- (40:2.7)
	(30:2.7) arc (30:150:2.7)
	(60:2.5) -- (64:2.5) -- (64:2.7)
	(40:2.5) -- (44:2.5) -- (44:2.7);

%% A-proj

\begin{scope}[xshift=-3cm]

\fill
	(40:2.7) circle (0.05);

\filldraw[fill=white]
	(1,{sqrt(3)}) circle (0.05);

\node at (-0.2,-0.1) {$O$};
\node at (0.85,1.95) {$A$};
\node at (42:2.95) {$B$};
\node at (60:2.95) {$M$};
\node at (138:2.9) {$B^*$};

\node at (0.63,1.5) {\scriptsize $\tfrac{1}{3}\pi$};
\node at (1.6,1.95) {\scriptsize $\tfrac{1}{n}\pi$};
\node at (0.2,1.95) {$D$};

\end{scope}

%% B-proj

\begin{scope}[xshift=3cm]

\filldraw[fill=white]
	(140:2.7) circle (0.05);

\fill
	(-1,{sqrt(3)}) circle (0.05);

\node at (0.2,-0.1) {$O$};
\node at (-0.85,1.95) {$B$};
\node at (41:2.98) {$A^*$};
\node at (120:2.95) {$M$};
\node at (138:2.95) {$A$};

\node at (-0.56,1.5) {\scriptsize $\tfrac{1}{n}\pi$};
\node at (-1.52,1.95) {\scriptsize $\tfrac{1}{3}\pi$};
\node at (0.2,1.95) {$D$};

\end{scope}
	
\end{tikzpicture}
\caption{Calculate $\triangle ABM$ in the $A$- and $B$-projections.}
\label{ABM}
\end{figure}

To get $d_{BM}$, we use the second of Figure \ref{ABM}, which is the triangle $\triangle ABM$ in the $B$-projection.  By
\[
A_B=-d_{AB},\;
A^*_B=d_{AB}^{-1},\;
M_B=d_{BM}e^{i(1-\frac{1}{n})\pi},
\]
we similarly get 
\begin{align*}
d_{BM}
&=|OM|_B-|OB|_B
=|OA|_B-|OB|_B \\
&=|BD|_B\sec(\tfrac{1}{2}-\tfrac{1}{3})\pi
-|AD|_B\sec\tfrac{1}{n}\pi  \\
&=\tfrac{1}{2}(d_{AB}^{-1}+d_{AB})\sec\tfrac{1}{6}\pi
-\tfrac{1}{2}(d_{AB}^{-1}-d_{AB})\sec\tfrac{1}{n}\pi. 
\end{align*}
Again the explicit values are given in Table \ref{distance}.

\subsection{Determination of Boundary}

In Figure \ref{boundary-tetraA}, for $v\in \Omega_5$, we connect $Av$ and rotate by $-\frac{2}{3}\pi$ around $A$ to get $Aw$. The picture is for the tetrahedron in $A$-projection, and is schematically correct for the octahedron and icosahedron. In this construction, we have $Av=a_1$, $Aw=a_2$ and $Bv=b_1$.

\begin{figure}[htp]
\centering
\begin{tikzpicture}[>=latex,scale=1.5]

\coordinate (B) at (1,0);
	
\foreach \a in {1,-1}
\draw[yscale=\a,gray] 
	(B) arc (30:160:{sqrt(3)})
	(B) arc (30:20:{sqrt(3)})
	(60:0.9) -- (240:1.2);

\draw[gray]
	(-2,0) -- (1,0);

\draw[dashed]
	(B) arc (60:120:3);

\draw
	(0,0) -- (139:0.61)
	(0,0) -- (139-120:0.61)
	(B) arc (60:90:3);

\filldraw[fill=white]
	(0,0) circle (0.035);

\fill
	(B) circle (0.035);

\node at (-0.03,-0.2) {$A$};
\node at (-50:0.88) {$C$};
\node at (1.15,0) {$B$};
\node at (-2.2,0) {$B^*$};

\node at (135:0.68) {$v$};
\node at (25:0.68) {$w$};

\node at (0.4,0.5) {\small $1$};
\node at (0,0.6) {$3$};
\node at (-1,0.6) {$5$};

\node at (0.4,-0.3) {$2$};
\node at (0,-0.55) {$4$};
\node at (-1,-0.5) {$6$};

\end{tikzpicture}
\caption{Boundary of the moduli space part in $\Omega_5$.}
\label{boundary-tetraA}
\end{figure}

If $Aw$ intersects $Bv$, then the pentagon is not simple, and $v$ is not in the moduli space. If $Aw$ does not intersect $Bv$, then by $b_2\in B-8-10$, $c_1\in C-2-1$, and $c_2\in C-10$, we know the pentagon is simple. The critical case between the two cases is when $B,w,v,B^*$ are on the same great circle. This means that the end $w$ of $a_2$ touches $b_1$. 

Figure \ref{curve} gives the general description. In the $X$-projection, let $P$ be the point at $d$ satisfying $0<d<1$. Then the antipode $P^*=-d^{-1}$. We also fix a rotation angle $2\alpha$, with $\alpha\in (0,\frac{1}{2}\pi)$. Then we look for $v=re^{i\phi}$, $w=ve^{-i2\alpha}=re^{i(\phi-2\alpha)}$, $r>0$, such that $P,w,v,P^*$ are on the same great circle. We assume that both $v$ and $w$ are in the upper half, which means $\phi\in [2\alpha,\pi]$. We also assume that the arc $PwvP^*$ is the ``shorter'' (actually the same length on the sphere) of the two arcs connecting $P$ and $P^*$. This corresponds to $\rho<\frac{1}{2}\pi$ in Lemma \ref{geometry1}. By the lemma, we must have $\phi\in [2\alpha,\tfrac{1}{2}\pi+\alpha]$.

\begin{figure}[htp]
\centering
\begin{tikzpicture}[>=latex,scale=2]

\node at (0,-0.15) {$X$};
\node at (1.2,-0.15) {$d$};
\node at (-1.8,-0.15) {$-d^{-1}$};
\node at (1.3,0.1) {$P$};
\node at (-1.9,0.1) {$P^*$};

\coordinate (P1) at (1.2,0);
\coordinate (P2) at (-1.8,0);
\coordinate (V) at (140:0.85);
\coordinate (W) at (140-120:0.85);
\coordinate (V1) at (140:1.2);
\coordinate (V2) at (140:0.6);
	
\arcThroughThreePoints[dashed]{V}{P1}{P2};
\arcThroughThreePoints[dashed]{V1}{P1}{P2};
\arcThroughThreePoints[dashed]{V2}{P1}{P2};
\arcThroughThreePoints{P1}{P2}{V};
\arcThroughThreePoints{P1}{P2}{V1};
\arcThroughThreePoints{P1}{P2}{V2};

\draw
	(-2,0) -- (1.5,0)
	(0,0) -- (140:1.2)
	(0,0) -- (20:1.01);

\filldraw[fill=white]
	(0,0) circle (0.035);

\node at (-0.64,0.61) {\small $v$};
\node at (0.8,0.37) {\small $w$};
\node at (-0.46,0.29) {\small $v_1$};
\node at (0.67,0.14) {\small $w_1$};
\node at (-0.95,0.86) {\small $v_2$};
\node at (1.05,0.37) {\small $w_2$};

\node at (0.05,0.15) {\small $2\alpha$};

\end{tikzpicture}
\caption{Determine the boundary curve of the moduli space.}
\label{curve}
\end{figure}

\begin{lemma}\label{geometry2}
In Figure \ref{curve}, suppose $XP<XP^*$, and $P,w,v,P^*$ are on the same great circle. Suppose $v_1,v_2$ are two points on the ray $Xv$, such that $v_1$ is within $v$ and $v_2$ is beyond $v$, and $\angle XPv_2<\frac{1}{2}\pi$. Suppose $w_1,w_2$ are the intersections of $Pv_1,Pv_2$ with the ray $Xw$. Then $Xv_1<Xw_1$ and $Xv_2>Xw_2$.
\end{lemma}

In the lemma, the inequalities in the Euclidean metric are equivalent to the same inequalities in the spherical metric. Therefore the lemma holds on the sphere.

The proof of the lemma is given by Figure \ref{geometry2_proof}. The first picture expands the spheres  $\bigcirc Pv_1P^*,\bigcirc PvP^*,\bigcirc Pv_2P^*$ in Figure \ref{curve}. Their respective centers $O_1,O,O_2$ are arranged on the same straight line as in the picture. Then we get the situation in the second picture, in which the circle is $\bigcirc Pv_1P^*$ or $\bigcirc Pv_2P^*$, and we may assume $OX$ is a horizontal line above the center of the circle. From $X$, we shoot two rays at angles $\pm\alpha$ with respect to $OX$, and the two rays intersect the circle at $V,W$. Then Lemma \ref{geometry2} effectively means $XV>XW$. To see the inequality, we flip the circle with respect to $OX$ and get a new circle. Since $OX$ is above the center of the circle, the new circle intersects $XV$ at $W'$. Then we get $XW=XW'<XV$. This completes the proof of Lemma \ref{geometry2}.

\begin{figure}[htp]
\centering
\begin{tikzpicture}[>=latex,scale=1]

\draw
	({sqrt(5)},0) -- (-{sqrt(5)},0)
	(-1,2) -- (-1,-2)
	(-1,{1/sqrt(3)}) -- ++(190:1.65)
	(-1,{1/sqrt(3)}) -- ++(110:1.35)
	;

\foreach \a in {0.9,-0.1,-0.7}
\draw
	(\a,0.1) -- ++(0.2,-0.2)
	(\a,-0.1) -- ++(0.2,0.2);
	
\draw[dashed]
	(0,0) circle ({sqrt(5)})
	(1,0) circle ({sqrt(8)})
	(-0.6,0) circle ({sqrt(4.16)});

\draw[densely dotted]
	(150:-3.66) -- (150:2.55);

\coordinate (P1) at (-1,2);
\coordinate (P2) at (-1,-2);
\coordinate (V) at (-2.21,0.36);
\coordinate (V1) at (-1.795,0.44);
\coordinate (V2) at (-2.62,0.295);

\arcThroughThreePoints{P1}{P2}{V};
\arcThroughThreePoints{P1}{P2}{V1};
\arcThroughThreePoints{P1}{P2}{V2};

\filldraw[fill=white]
	(-1,{1/sqrt(3)}) circle (0.05);

\node at (0,-0.3) {\small $O$};
\node at (1.1,-0.3) {\small $O_1$};
\node at (-0.55,-0.3) {\small $O_2$};
\node at (-0.8,0.55) {\small $X$};

\node at (-1.1,2.2) {\small $P$};
\node at (-1.15,-2.2) {\small $P^*$};

\node at (-1.4,0.65) {\small $\alpha$};
\node at (-1.3,0.9) {\small $\alpha$};

\node at (-2.33,0.47) {\small $v$};
\node at (-1.6,0.3) {\small $v_1$};
\node at (-2.85,0.3) {\small $v_2$};
\node at (-2.2,1.7) {\small $w$};

\draw[->]
	(-2,1.7) -- (-1.42,1.73);

\begin{scope}[xshift=7cm]

\draw[dashed]
	(0,0) circle (2)
	(2,1.6) arc (0:-180:2);

\draw[densely dotted]
	(-1.833,0.8) -- (1.833,0.8);

\draw
	(50:2) -- (50:-2)
	(-0.5,0.8) -- ++(140:1.08)
	(-0.5,0.8) -- ++(220:1.9)
	(0.67,0.8) -- ++(95:0.14)
	(0.67,0.8) -- ++(5:0.14)
	(0.67,0.8) -- ++(95:-0.14)
	(0.67,0.8) -- ++(5:-0.14)
	(0,0) -- ++(95:0.14)
	(0,0) -- ++(5:0.14)
	(0,0) -- ++(95:-0.14)
	(0,0) -- ++(5:-0.14);

\filldraw[fill=white]
	(-0.5,0.8) circle (0.05);
	
\node at (-0.9,0.95) {\small $\alpha$};
\node at (-0.9,0.65) {\small $\alpha$};
	
\node at (0.85,0.6) {\small $O$};
\node at (-0.35,0.6) {\small $X$};
\node at (130:2.2) {\small $W$};
\node at (190:2.2) {\small $V$};
\node at (-1.2,-0.2) {\small $W'$};
\node at (0.8,-0.2) {\small $O_1$ or $O_2$};

\end{scope}

\end{tikzpicture}
\caption{Proof of Lemma \ref{geometry2}.}
\label{geometry2_proof}
\end{figure}

We apply Lemma \ref{geometry2} to Figure \ref{boundary-tetraA}. For any $\phi\in [2\alpha,\tfrac{1}{2}\pi+\alpha]$, $\alpha=\frac{1}{3}\pi$, we consider the ray from $A$ at angle $\phi=\angle BAv$. There is a point $v$ on the ray, such that $B,w,v,B^*$ are on the same circle. Then Lemma \ref{geometry2} says that a point on the ray is inside the moduli space if and only if it is between $A$ and $v$. Therefore the track of $v$ in Figure \ref{boundary-tetraA} gives a curve $\gamma_A$. The part of the moduli space in $\Omega_5$ is the region between $A$ and $\gamma_A$ (not including the curve), and is bounded by $\gamma_A$ and the arc $\Omega_3\cap\Omega_5$ (including the interior of the arc).

\subsection{Formula for Boundary}

In Figure \ref{curve}, the four points $v=re^{i\phi}$, $w=ve^{-i2\alpha}=re^{i(\phi-2\alpha)}$, $P=d$, $P^*=-d^{-1}$ are on the same great circle if and only if the following is a real number ($=_{\bb R}$ means equality up to multiplying and adding real numbers)
\begin{align*}
\frac{v+d^{-1}}{v-w}\cdot \frac{d-w}{d+d^{-1}}
&=_{\bb R} \frac{(v+d^{-1})(d-ve^{-i2\alpha})}{v(1-e^{-i2\alpha})} \\
&=_{\bb R} i((d-d^{-1})\cos\alpha-ve^{-i\alpha}+v^{-1}e^{i\alpha}).
\end{align*}
The condition for the number to be real is a cubic equation
\begin{equation}\label{eqA}
4\lambda |v|^2+(ve^{-i\alpha}+\bar{v}e^{i\alpha})(|v|^2-1)=0,
\end{equation}
where
\begin{equation}\label{eqF}
\lambda=\tfrac{1}{2}(d^{-1}-d)\cos\alpha>0.
\end{equation}
In terms of $\xi=(\xi_1,\xi_2,\xi_3)$ on the sphere, the equation \eqref{eqA} for the track $\gamma$ of $v$ becomes quadratic
\begin{equation}\label{eqB}
\lambda(1-\xi_3^2)+(\xi_1\cos\alpha+\xi_2\sin\alpha)\xi_3
=\lambda(\xi_1^2+\xi_2^2)+(\xi_1\cos\alpha+\xi_2\sin\alpha)\xi_3
=0.
\end{equation}
The first version is a hyperbolic cylinder, and the second version is homogeneous. Therefore the curve is the intersection of a hyperbolic cylinder with the sphere, and is also the intersection of a quadratic cone with the sphere.

For $v=x+iy$, \eqref{eqA} becomes the cartesian form of $\gamma$
\begin{equation}\label{eqC}
2\lambda(x^2+y^2)+(x\cos\alpha+y\sin\alpha)(x^2+y^2-1)=0.
\end{equation}
For $v=re^{i\phi}$, $\gamma$ is a curve connecting $de^{i2\alpha}$ to $0$, given by (see \eqref{eqR1})
\begin{equation}\label{eqE}
r^{-1}-r
=2R=2\lambda\sec(\phi-\alpha),\quad
\phi\in [2\alpha,\tfrac{1}{2}\pi+\alpha].
\end{equation}
The formula implies that $R>0$ and is strictly increasing in $\phi$. This implies that $1>r>0$, and $r$ is strictly decreasing in $\phi$. Then we take the non-negative solution of \eqref{eqE} to get the explicit formula
\begin{equation}\label{eqD}
r=\sqrt{\lambda^2\sec^2(\phi-\alpha)+1}-\lambda\sec(\phi-\alpha),\quad
\phi\in [2\alpha,\tfrac{1}{2}\pi+\alpha].
\end{equation}
By \eqref{eqR2}, we also get the spherical version of the curve
\begin{equation}\label{eqG}
\xi=\frac{(\cos\phi,\sin\phi,-\lambda\sec(\phi-\alpha))}{\sqrt{\lambda^2\sec^2(\phi-\alpha)+1}},\quad
\phi\in [2\alpha,\tfrac{1}{2}\pi+\alpha].
\end{equation}
This is equivalent to \eqref{eqB}.

We apply the calculation to $\gamma_A$. This means $X=A$, $P=B$, $\alpha=\frac{1}{3}\pi$, $d=d_{AB}$. Moreover, $z_A=re^{i\theta}$ equals $v=re^{i\phi}$. Substituting the data into \eqref{eqC} through \eqref{eqG}, we get $\lambda_A$ in Table \ref{coefficient} and the complex equation for the $A$-projection of $\gamma_A$
\[
4\lambda_A|z|^2+(ze^{-i\frac{1}{3}\pi}+\bar{z}e^{i\frac{1}{3}\pi})(|z|^2-1)=0.
\]
The cartesian equation for $z=x+iy$ is
\[
4\lambda_A(x^2+y^2)+(x+\sqrt{3}y)(x^2+y^2-1)=0.
\]
The polar equation for $z=re^{i\theta}$ is
\[
r=\sqrt{\lambda_A^2\sec^2(\theta-\tfrac{1}{3}\pi)+1}-\lambda_A\sec(\theta-\tfrac{1}{3}\pi),\quad
\theta\in [\tfrac{2}{3}\pi,\tfrac{5}{6}\pi].
\]
The spherical equation is
\[
\xi=\frac{(\cos\theta,\sin\theta,-\lambda_A\sec(\theta-\tfrac{1}{3}\pi))}{\sqrt{\lambda_A^2\sec^2(\theta-\tfrac{1}{3}\pi)+1}},\quad
\theta\in [\tfrac{2}{3}\pi,\tfrac{5}{6}\pi].
\]

\begin{table}[htp]
\centering
\begin{tabular}{|l||c|c|c|}
\hline 
& $n=3$   
& $n=4$  
& $n=5$  \\
\hline \hline
$\lambda_A=\tfrac{1}{2}(d_{AB}^{-1}-d_{AB})\cos\tfrac{1}{3}\pi$
& $\frac{1}{4\sqrt{2}}$
& $\frac{1}{2\sqrt{2}}$ 
& $\frac{\sqrt{5}+3}{8}$ \\
\hline 
$\lambda_B=\tfrac{1}{2}(d_{AB}^{-1}-d_{AB})\cos\tfrac{1}{n}\pi$  
& $\frac{1}{4\sqrt{2}}$
& $\frac{1}{2}$ 
& $\frac{\sqrt{5}+2}{4}$ \\
\hline 
$\lambda_C=\tfrac{1}{2}(d_{AM}^{-1}-d_{AM})\cos\tfrac{1}{3}\pi$  
& $\frac{1}{2\sqrt{2}}$
& $\frac{1}{\sqrt{2}}$ 
& $\frac{\sqrt{5}+3}{4}$  \\
\hline  
\end{tabular}
\caption{Coefficients for $\gamma_A$, $\gamma_B$, $\gamma_C$.}
\label{coefficient}
\end{table}

Next we apply the calculation to $\Omega_{13}$ (for $n=3,4$) or $\Omega_{13}\cup\Omega_{19}$ (for $n=5$). The situation is similar to $\Omega_5$, with $a$ and $b$ exchanged. This means that $A$ and $B$ are exchanged, and the calculation can be done in the $B$-projection. We denote the track by $\gamma_B$, and $V\in \Omega_{13}$ (or $\Omega_{13}\cup\Omega_{19}$) gives a simple pentagon if and only if $V$ lies between $\gamma_B$ (not including the curve) and the arc $\Omega_7\cap\Omega_{13}$.  

By comparing the $B$-projection with the $A$-projection, we may calculate $\gamma_B$ by horizontally flipping the standard picture in Figure \ref{curve}. This means taking $X=B$, $P=A$, $\alpha=\frac{1}{n}\pi$, $d=d_{AB}$. Moreover, for $z_B=re^{i\theta}$, we take $v=re^{i(\pi-\theta)}=-\bar{z}_B$. Substituting the data into \eqref{eqC} through \eqref{eqG}, we get $\lambda_B$ in Table \ref{coefficient} and the complex equation for the $B$-projection of $\gamma_B$
\[
4\lambda_B|z|^2-(ze^{i\frac{1}{n}\pi}+\bar{z}e^{-i\frac{1}{n}\pi})(|z|^2-1)=0.
\]
The cartesian equation is
\[
2\lambda_B(x^2+y^2)-(x\cos\tfrac{1}{n}\pi-y\sin\tfrac{1}{n}\pi)(x^2+y^2-1)=0.
\]
The polar equation is
\[
r=\sqrt{\lambda_B^2\sec^2(\theta+\tfrac{1}{n}\pi)+1}+\lambda_B\sec(\theta+\tfrac{1}{n}\pi),\quad
\theta\in [(\tfrac{1}{2}-\tfrac{1}{n})\pi,(1-\tfrac{2}{n})\pi].
\]
The spherical equation is
\[
\xi=\frac{(\cos\theta,\sin\theta,\lambda_B\sec(\theta+\tfrac{1}{n}\pi))}{\sqrt{\lambda_B^2\sec^2(\theta+\tfrac{1}{n}\pi)+1}},\quad
\theta\in [(\tfrac{1}{2}-\tfrac{1}{n})\pi,(1-\tfrac{2}{n})\pi].
\]

Next we apply the calculation to $\Omega_4$. This is described by Figure \ref{boundary-tetraC}, which is the rotation of the first of Figure \ref{projection-tetra} by $-\frac{2}{3}\pi$. For $v\in \Omega_4$, we rotate $v$ by $-\frac{2}{3}\pi$ around $A$ to get $w$, and require $C,v,w,C^*$ to be on the same great circle. If $V=v$, then $Av=a_1, Aw=a_2,Cw=c_1$, and the picture describes the extreme case that the end $v$ of $a_1$ touches $c_1$. The track of $v$ gives a curve $\gamma_C$. By Lemma \ref{geometry2}, a point $V\in \Omega_4$ gives a simple pentagon if and only if $V$ lies between the curve $\gamma_C$ (not including the curve) and the arc $AC$.

\begin{figure}[htp]
\centering
\begin{tikzpicture}[>=latex,scale=2]

\coordinate (C) at ({(-sqrt(3)+1)/sqrt(2)},0);
\coordinate (D) at ({(sqrt(3)+1)/sqrt(2)},0);
\coordinate (X) at (240:{1/sqrt(2)});
\pgfmathsetmacro{\r}{sqrt(3)/sqrt(2)};	

\draw[gray] 
	(X) arc (-90:5:{\r})
	(X) arc (-90:-100:{\r})
	(X) arc (210:145:{\r})
	(X) arc (210:215:{\r})
	(D) arc (0:35:{\r})
	(D) arc (0:-35:{\r})
	(60:0.8) -- (240:0.8)
	(120:0.8) -- (-60:0.8)
	(-0.8,0) -- (2.2,0);

\draw[dashed]
	(C) arc (115:65:2.9);
	
\draw
	(C) arc (115:99.5:2.9)
	(0,0) -- (165:0.33)
	(0,0) -- (165-120:0.33);

\filldraw[fill=white]
	(0,0) circle (0.025);

\fill
	(X) circle (0.025);

\node at (2.05,0.1) {$C^*$};
\node at (-0.6,0.1) {$C$};
\node at (-0.03,-0.15) {$A$};
\node at (-0.2,-0.5) {$B$};

\node at (-0.32,0.15) {\small $v$};
\node at (0.27,0.3) {\small $w$};

\node at (0,-0.4) {$1$};
\node at (0.35,-0.2) {$3$};
\node at (0.6,0.5) {$5$};
\node at (0,0.6) {$6$};
\node at (-0.35,0.35) {$4$};
\node at (-0.35,-0.25) {$2$};

\end{tikzpicture}
\caption{Deriving the boundary of the moduli space part in $\Omega_4$.}
\label{boundary-tetraC}
\end{figure}

Matching Figure \ref{boundary-tetraC} with Figure \ref{curve} means taking $X=A$, $P=C$, $\alpha=\frac{1}{3}\pi$, $d=d_{AC}=d_{AM}$. Moreover, for $z_A=re^{i\theta}$, we take $v=re^{i(\frac{1}{3}\pi-\theta)}=\bar{z}_Ae^{i\frac{1}{3}\pi}$ (note that $v$ and $w$ are exchanged in the matching of the two pictures). Substituting the data into \eqref{eqC} through \eqref{eqG}, we get $\lambda_C$ in Table \ref{coefficient} and the complex equation for the $A$-projection of $\gamma_C$
\[
4\lambda_C|z|^2+(z+\bar{z})(|z|^2-1)=0.
\]
The cartesian equation is
\[
2\lambda_C(x^2+y^2)+x(x^2+y^2-1)=0.
\]
The polar equation is
\[
r=\sqrt{\lambda_C^2\sec^2\theta+1}-\lambda_C\sec\theta,\quad
\theta\in [-\tfrac{1}{2}\pi,0].
\]
The spherical equation is
\[
\xi=\frac{(\cos\theta,\sin\theta,-\lambda_C\sec\theta)}{\sqrt{\lambda_C^2\sec^2\theta+1}},\quad
\theta\in [-\tfrac{1}{2}\pi,0].
\]

We may further calculate the boundary of the moduli space part in $\Omega_8$ ($n=3,4$) or $\Omega_8\cup\Omega_{14}$ ($n=5$), in the $B$-projection. The coefficient
\[
\lambda
=\tfrac{1}{2}(d_{BM}^{-1}-d_{BM})\cos\tfrac{1}{n}\pi
=\tfrac{1}{2}(d_{AM}^{-1}-d_{AM})\cos\tfrac{1}{3}\pi
=\lambda_C.
\]
The complex equation for the $B$-projection of $\gamma_C$ is
\[
4\lambda_C|z|^2-(z+\bar{z})(|z|^2-1)=0,
\]
The cartesian equation is
\[
2\lambda_C(x^2+y^2)-x(x^2+y^2-1)=0,
\]
The polar equation is
\[
r=\sqrt{{\lambda_C}^2\sec^2\theta+1}+\lambda_C\sec\theta,\quad
\theta\in [-\pi,-\tfrac{1}{2}\pi].
\]
The spherical equation is
\[
\xi=\frac{(\cos\theta,\sin\theta,\lambda_C\sec\theta)}{\sqrt{\lambda_C^2\sec^2\theta+1}},\quad
\theta\in [-\pi,-\tfrac{1}{2}\pi].
\]

We note that the $A$-projection and the $B$-projection for $\gamma_C$ are related by the transformation $(x,y)\mapsto (-x,y)$. Geometrically, this means that $\gamma_C$ is symmetric with respect to the great circle of equal distance from $A$ and $B$. 

We use the same notation $\gamma_C$ for the boundary curve in $\Omega_4$ and in $\Omega_8$ (or $\Omega_8\cup\Omega_{14}$), because the two curves are actually the same one. Figure \ref{same_curve} gives the schematic reason for this claim. The curves in the two regions describe the extreme case that $v$ lies on $c$ ($c_1$ for $\Omega_4$ and $c_2$ for $\Omega_8$ (or $\Omega_8\cup\Omega_{14}$)). Although $w$ is an end of $a_2$ for the first region and an end of $b_2$ for the second region, the two versions of $w$ and $C,C^*$ are on the same great circle. Therefore the formulae for the curves are obtained by the same requirement that $C,v,w,C^*$ are on the same great circle, irrespective to which $w$ we choose. 

In the polar equation for the $A$-projection of $\gamma_C$, the curve is in $\Omega_4$ for $\theta\in [-\tfrac{1}{2}\pi,-\tfrac{1}{3}\pi]$, and in $\Omega_8$ (or $\Omega_8\cup\Omega_{14}$) for $\theta\in [-\tfrac{1}{3}\pi,0]$. The two combine to yield the range $[-\tfrac{1}{2}\pi,0]$ we gave earlier. Similar remark can be made for the polar equation for the $B$-projection of $\gamma_C$.

\begin{figure}[htp]
\centering
\begin{tikzpicture}[>=latex,scale=1]

\draw
	(-3,0) -- (-1.5,1) -- (0,0)
	(-0.5,-0.1) -- (-0.5,0.1);

\draw[dashed]
	(-3,0) -- (2,0);

\draw[line width=1.2]
	(2,0) -- (1,1.5) -- (0,0);

\draw[->]
	(2.2,0) -- ++(0.4,0);
	
\node at (-2.3,0.75) {$a_2$};
\node at (-0.7,0.75) {$a_1$};
\node at (0.3,0.9) {$b_1$};
\node at (1.7,0.9) {$b_2$};
\node at (-1.5,-0.2) {$c_1$};
\node at (1,-0.2) {$c_2$};

\node at (-1.5,1.2) {$A$};
\node at (1,1.7) {$B$};
\node at (-0.5,-0.3) {$C$};
\node at (3,0) {$C^*$};
\node at (0,-0.2) {$v$};
\node at (-3,-0.2) {$w$};
\node at (2,-0.2) {$w$};

\end{tikzpicture}
\caption{Moduli space parts in $\Omega_4$ and $\Omega_8$ (or $\Omega_8\cup\Omega_{14}$) have the same boundry.}
\label{same_curve}
\end{figure}

\section{Picture of Moduli Space}
\label{picture}

By the discussion in Sections \ref{region} and \ref{boundary}, the moduli space is the open subset of the sphere bounded by 
\begin{itemize}
\item arcs $\Omega_3\cap\Omega_9$ and $\Omega_7\cap\Omega_9$;
\item curve $\gamma_A$ connecting two ends of $\Omega_3\cap\Omega_5$, and inside $\Omega_5$;
\item curve $\gamma_B$ connecting two ends of $\Omega_7\cap\Omega_{13}$, and inside $\Omega_{13}$ (or $\Omega_{13}\cup\Omega_{19}$);
\item curve $\gamma_C$ connecting $A$ to $C$ to $B$, and inside $\Omega_4$ and $\Omega_8$ (or $\Omega_8\cup\Omega_{14}$). 
\end{itemize}
Figure \ref{3curves} gives the perspective pictures of the three curves in the most convenient projections. Recall that the numbers 3, 4, 5 mean tetrahedron, octahedron, icosahedron. Moreover, the point $B'$ is obtained by rotating $B$ around $A$ by $\frac{2}{3}\pi$, and $A'$ is obtained by rotating $A$ around $B$ by $-\frac{2}{n}\pi$. 

\begin{figure}[htp]
\centering
\begin{tikzpicture}[>=latex,scale=3]

%% \gamma_A

\draw[gray]
	(-0.5,0) -- (0.8,0)
	(0,0) -- (120:0.8);

\draw
	plot [domain=60:89.5, samples=100] ({60+\x}:{sqrt(1+(1/32)*sec(\x)*sec(\x))-(1/(4*sqrt(2)))*sec(\x)})
	plot [domain=60:89.5, samples=100] ({60+\x}:{sqrt(1+(1/8)*sec(\x)*sec(\x))-(1/(2*sqrt(2)))*sec(\x)})
	plot [domain=60:89.5, samples=100] ({60+\x}:{sqrt(1+((sqrt(5)+3)/8)*((sqrt(5)+3)/8)*sec(\x)*sec(\x))-((sqrt(5)+3)/8)*sec(\x)})
	plot [domain=-89.5:0, samples=100] ({\x}:{sqrt(1+(1/8)*sec(\x)*sec(\x))-(1/(2*sqrt(2)))*sec(\x)})
	plot [domain=-89.5:0, samples=100] ({\x}:{sqrt(1+(1/2)*sec(\x)*sec(\x))-(1/sqrt(2))*sec(\x)})
	plot [domain=-89:0, samples=100] ({\x}:{sqrt(1+1.7135*sec(\x)*sec(\x))-1.3090*sec(\x)});

\foreach \a in {0,120}
\fill
	(\a:{1/sqrt(2)}) circle (0.02)
	(\a:{(sqrt(3)-1)/sqrt(2)}) circle (0.02)
	(\a:0.3382612127) circle (0.02);
	
\filldraw[fill=white]
	(0,0) circle (0.02);

\node at (120:0.88) {\scriptsize $\frac{2}{3}\pi$};
	
\node at (0.08,0.07) {$A$};	
\node at (0.5,0.08) {$B$};	
\node at (-0.15,0.5) {$B'$};

\node at (-0.38,0.2) {$\gamma_A$};	
\node at (0,-0.35) {$\gamma_C$};	

\node at (0.65,-0.35) {\small $3$};
\node at (0.47,-0.25) {\small $4$};
\node at (0.35,-0.15) {\small $5$};

%% \gamma_B

\begin{scope}[xshift=2cm]

\draw[gray]
	(-0.8,0) -- (0.5,0)
	(0,0) -- (60:0.8)
	(0,0) -- (90:0.8)
	(0,0) -- (108:0.8);	
	
\draw
	plot [domain=90.5:120, samples=100] ({-60+\x}:{sqrt(1+(1/32)*sec(\x)*sec(\x))+(1/(4*sqrt(2)))*sec(\x)})
	plot [domain=90.5:135, samples=100] ({-45+\x}:{sqrt(1+(1/4)*sec(\x)*sec(\x))+(1/2)*sec(\x)})
	plot [domain=90.5:144, samples=100] ({-36+\x}:{sqrt(1+1.1215169944*sec(\x)*sec(\x))+1.0590169944*sec(\x)})
	plot [domain=-90.5:-180, samples=100] ({\x}:{sqrt(1+(1/8)*sec(\x)*sec(\x))+(1/(2*sqrt(2)))*sec(\x)})
	plot [domain=-90.5:-180, samples=100] ({\x}:{sqrt(1+(1/2)*sec(\x)*sec(\x))+(1/sqrt(2))*sec(\x)})
	plot [domain=-91:-180, samples=100] ({\x}:{sqrt(1+1.7135*sec(\x)*sec(\x))+1.3090*sec(\x)});

\filldraw[fill=white]
	(60:{1/sqrt(2)}) circle (0.02)
	(90:{(sqrt(3)-1)/sqrt(2)}) circle (0.02)
	(108:0.3382612127) circle (0.02)
	(180:{1/sqrt(2)}) circle (0.02)
	(180:{(sqrt(3)-1)/sqrt(2)}) circle (0.02)
	(180:0.3382612127) circle (0.02);
	
\fill
	(0,0) circle (0.02);

\node at (-0.1,0.07) {$B$};	
\node at (0.1,0.6) {$A'$};	
\node at (-0.5,0.08) {$A$};

\node at (0.4,0.2) {$\gamma_B$};	
\node at (0,-0.35) {$\gamma_C$};	

\node at (-0.65,-0.35) {\small $3$};
\node at (-0.47,-0.25) {\small $4$};
\node at (-0.35,-0.15) {\small $5$};

\node at (60:0.9) {\scriptsize $\frac{1}{3}\pi$};
\node at (90:0.85) {\scriptsize $\frac{1}{2}\pi$};
\node at (108:0.85) {\scriptsize $\frac{3}{5}\pi$};

\end{scope}
	
\end{tikzpicture}
\caption{$A$-projection of $\gamma_A,\gamma_C$, and $B$-projection of $\gamma_B,\gamma_C$.}
\label{3curves}
\end{figure}

\begin{figure}[htp]
\centering
\begin{tikzpicture}[>=latex,scale=5]

\draw[samples=100, smooth, domain=90:180]
	plot ({\x}:{( % A3
		sqrt(
			(sqrt(6*cos(\x)*cos(\x)+2)+sqrt(2)*cos(\x))*
			(sqrt(6*cos(\x)*cos(\x)+2)+sqrt(2)*cos(\x))+4
			) 
		- sqrt(6*cos(\x)*cos(\x)+2) - sqrt(2)*cos(\x)
			)/2})
	plot ({\x}:{sqrt(cos(\x)*cos(\x)+2)-sqrt(cos(\x)*cos(\x)+1)}) % A4
	plot ({\x}:{( % A5
		sqrt( 
			(sqrt((10-2*sqrt(5))*cos(\x)*cos(\x)+2*sqrt(5)+6)
				-(sqrt(5)-1)*cos(\x))*
			(sqrt((10-2*sqrt(5))*cos(\x)*cos(\x)+2*sqrt(5)+6)
				-(sqrt(5)-1)*cos(\x)) + 4 
			)  
		-(sqrt((10-2*sqrt(5))*cos(\x)*cos(\x)+2*sqrt(5)+6)
				-(sqrt(5)-1)*cos(\x)) 
			)/2});

\draw[samples=100, smooth, domain=0:-90]
	plot ({\x}:{( % B3
		sqrt(
			(sqrt(6*sin(\x)*sin(\x)+2)+sqrt(2)*sin(\x))*
			(sqrt(6*sin(\x)*sin(\x)+2)+sqrt(2)*sin(\x))+4
			) 
		- sqrt(6*sin(\x)*sin(\x)+2) - sqrt(2)*sin(\x)
			)/2})
	plot ({\x}:{ % B4
		sqrt(
			(sqrt(2*sin(\x)*sin(\x)+2)+sin(\x))*
			(sqrt(2*sin(\x)*sin(\x)+2)+sin(\x))+1
			) 
		- sqrt(2*sin(\x)*sin(\x)+2) - sin(\x)
			})
	plot ({\x}:{( % B5
		sqrt( 
			(sqrt((10+2*sqrt(5))*sin(\x)*sin(\x)+6*sqrt(5)+14)
				+(sqrt(5)+1)*sin(\x))*
			(sqrt((10+2*sqrt(5))*sin(\x)*sin(\x)+6*sqrt(5)+14)
				+(sqrt(5)+1)*sin(\x)) + 4 
			)  
		-(sqrt((10+2*sqrt(5))*sin(\x)*sin(\x)+6*sqrt(5)+14)
				+(sqrt(5)+1)*sin(\x)) 
			)/2});

\draw[samples=100, smooth, domain=180:270]
	plot ({\x}:{( % C3
		sqrt(
			(sqrt(2)*(2*cos(\x)*sin(\x)-1)/(cos(\x)+sin(\x)))*
			(sqrt(2)*(2*cos(\x)*sin(\x)-1)/(cos(\x)+sin(\x)))+4
			) 
		- (sqrt(2)*(2*cos(\x)*sin(\x)-1)/(cos(\x)+sin(\x)))
			)/2})
	plot ({\x}:{ % C4
		sqrt(
			((cos(\x)*sin(\x)-sqrt(2))/(cos(\x)+sqrt(2)*sin(\x)))*
			((cos(\x)*sin(\x)-sqrt(2))/(cos(\x)+sqrt(2)*sin(\x)))+1
			) 
		- (cos(\x)*sin(\x)-sqrt(2))/(cos(\x)+sqrt(2)*sin(\x)) 
			})
	plot ({\x}:{ % C5
		sqrt(
			((sqrt(5)-1)*(cos(\x)*sin(\x)-sqrt(5)-2)/(2*cos(\x)+(sqrt(5)+1)*sin(\x)))*
			((sqrt(5)-1)*(cos(\x)*sin(\x)-sqrt(5)-2)/(2*cos(\x)+(sqrt(5)+1)*sin(\x))) + 1
			) 
		- ((sqrt(5)-1)*(cos(\x)*sin(\x)-sqrt(5)-2)/(2*cos(\x)+(sqrt(5)+1)*sin(\x)))
			});

\coordinate (A31) at ({-(sqrt(3)-1)/(sqrt(2))},0);
\coordinate (B31) at (0,{-(sqrt(3)-1)/(sqrt(2))});
\coordinate (A41) at ({-sqrt(3)+sqrt(2)},0);
\coordinate (B41) at (0,{-sqrt(2)+1});
\coordinate (A51) at ({(-sqrt(15)+sqrt(5)-sqrt(3)+3)/2},0);
\coordinate (B51) at (0,{-(sqrt(10+2*sqrt(5))-sqrt(5)-1)/2});

\coordinate (A32) at ({(sqrt(3)-1)/(sqrt(2))},0);
\coordinate (B32) at (0,{(sqrt(3)-1)/(sqrt(2))});
\coordinate (A42) at ({sqrt(3)-sqrt(2)},0);
\coordinate (B42) at (0,{sqrt(2)-1});
\coordinate (A52) at ({-(-sqrt(15)+sqrt(5)-sqrt(3)+3)/2},0);
\coordinate (B52) at (0,{(sqrt(10+2*sqrt(5))-sqrt(5)-1)/2});

\draw
	(A32) -- (0,0) -- (B32);

\draw[gray!40]
	(A31) -- (0,0) -- (B31)
	(A31) arc (210:240:{sqrt(2)})
	(A32) arc (-30:-90:{sqrt(2)})
	(B32) arc (120:180:{sqrt(2)})
	(A41) arc (210:225:2)
	(A42) arc (-30:-45:2)
	(B42) arc (135:160.56:2)
	(B41) arc (-90:-109.52:{sqrt(2)})
	(A52) arc (-30:-36:{sqrt(5)+1})
	(A51) arc (210:216:{sqrt(5)+1})
	(B52) arc (144:153.45:{sqrt(5)+1})
	(B51) arc (252:243.4:2); 
		
\foreach \a in {3,4,5}
\foreach \b in {1,2}
{
\filldraw[fill=white]
	(A\a\b) circle (0.01);

\fill
	(B\a\b) circle (0.01);
}

\node at (-0.67,-0.67) {$3$};
\node at (-0.38,-0.4) {$4$};
\node at (-0.2,-0.24) {$5$};

\node at (-0.07,-0.1) {$\Omega_1$};
\node at (-0.07,0.1) {$\Omega_3$};
\node at (0.07,-0.1) {$\Omega_7$};

\node at (0.05,0.05) {$M$};
\node at (-0.57,0.03) {$A$};
\node at (0.58,0.03) {$A'$};
\node at (0.01,-0.57) {$B$};
\node at (0.03,0.57) {$B'$};

\node at (-0.4,0.48) {$\gamma_A$};
\node at (0.48,-0.4) {$\gamma_B$};
\node at (-0.74,-0.74) {$\gamma_C$};

\node at (1.2,0) 
	{\includegraphics[scale=0.2]{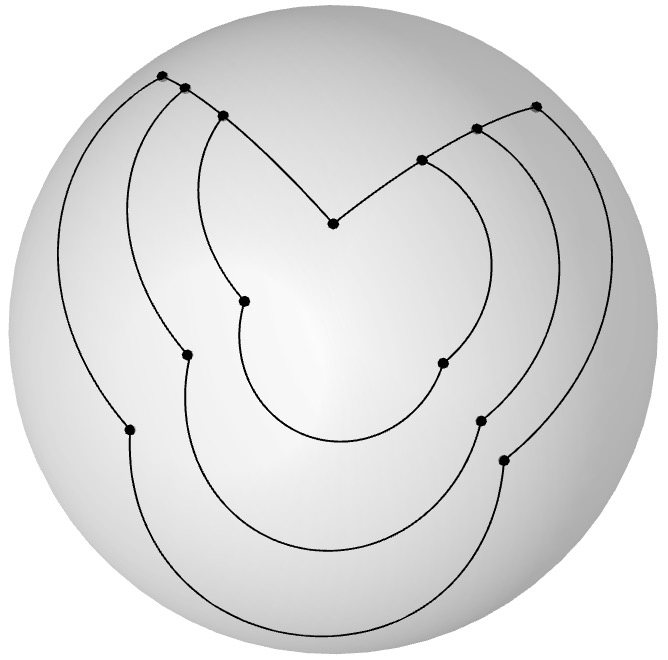}};

\end{tikzpicture}
\caption{Moduli space of pentagonal subdivision in $M$-projection, and in sphere.}
\label{moduliM}
\end{figure}

The most visually pleasing perspective picture of the moduli space is given by the $M$-projection in Figure \ref{moduliM}. We also indicate the four congruent triangular regions $\Omega_1,\Omega_2,\Omega_3,\Omega_7$ contained in the moduli space. In the $M$-projection, $\Omega_1,\Omega_3,\Omega_7$ are visually congruent, and $\Omega_2$ is the triangle outside $\Omega_1$. The following are some key points in the $M$-projection
\[
A_M=-d_{AM},\;
B_M=-d_{BM}i,\;
A'_M=d_{AM},\;
B'_M=d_{BM}i.
\]
To get the $M$-projection formulae for the three curves, we use the M\"obius transform. The M\"obius transform from the $M$-projection to the $A$-projection is easily determined by applying the transform to $A,A^*,M$ 
\[
A_M\to A_A=0,\quad 
A_M^*\to A^*_A=\infty,\quad
M_M=0\to M_A.
\]
We get
\begin{equation}\label{m2a}
z_A=M_A\frac{A^*_M}{A_M}\cdot\frac{z_M-A_M}{z_M-A^*_M}
=\frac{z_M+d_{AM}}{1-d_{AM}z_M}e^{i\frac{1}{3}\pi}.
\end{equation}
Substituting into the $A$-projection formulae of $\gamma_A$ and $\gamma_C$, we get the $M$-projection formulae of $\gamma_A$ and $\gamma_C$. Similarly, we get the $M$-projection formula of $\gamma_B$ by 
\begin{equation}\label{m2b}
z_B
=\frac{z_M+d_{BM}i}{1+d_{BM}iz_M}e^{i(\frac{1}{2}-\frac{1}{n})\pi}.
\end{equation}
The first of Figure \ref{moduliM} gives the perspective pictures of the $M$-projections of $\gamma_A,\gamma_B,\gamma_C$.

%%% moduli

The following are the cartesian forms of the three curves in the $M$-projection. For tetrahedron
\begin{align*}
\gamma_A &\colon
(x^2+y^2)^2+2\sqrt{2}(x^2+y^2)x-8x^2-4y^2-x+1=0, \\
\gamma_B &\colon
(x^2+y^2)^2+2\sqrt{2}(x^2+y^2)y-4x^2-8y^2-y+1=0, \\
\gamma_C &\colon
(x^2+y^2-1)(x+y)-\sqrt{2}(x-y)^2=0.
\end{align*}
For octahedron
\begin{align*}
\gamma_A &\colon
(x^2+y^2)^2-10x^2-6y^2+1=0, \\
\gamma_B &\colon
(x^2+y^2)^2+4(x^2+y^2)y-10x^2-14y^2-4y+1=0, \\
\gamma_C &\colon
(x^2+y^2-1)(x+\sqrt{2}y)-2\sqrt{2}(x^2+y^2)+2xy=0.
\end{align*}
For icosahedron
\begin{align*}
\gamma_A &\colon
(x^2+y^2)^2-2(\sqrt{5}-1)(x^2+y^2)x \\
&\quad -2(\sqrt{5}+6)x^2-2(\sqrt{5}+4)y^2+2(\sqrt{5}-1)x+1=0, \\
\gamma_B &\colon
(x^2+y^2)^2+2(\sqrt{5}+1)(x^2+y^2)y \\
&\quad -2(3\sqrt{5}+8)x^2-2(3\sqrt{5}+10)y^2-2(\sqrt{5}+1)y+1=0, \\
\gamma_C &\colon
(x^2+y^2-1)(2x+(\sqrt{5}+1)y) \\
&\quad -2(\sqrt{5}+3)(x^2+y^2)+2(\sqrt{5}-1)xy=0.
\end{align*}

Using \eqref{eqR1} and \eqref{eqR2}, the polar forms of the $M$-projections of the curves are given by ($R>0$ and $1>r>0$)
\[
r^{-1}-r=2R,\quad
r=\sqrt{R^2+1}-R,
\]
and the spherical forms of the curves are given by
\[
\xi=\frac{(\cos\theta,\sin\theta,-R)}{\sqrt{R^2+1}}.
\]
The function $R$ for the $M$-projection is explicitly given below 
\begin{align*}
R_A
&=\begin{cases}
\tfrac{1}{\sqrt{2}}\cos\theta+\tfrac{1}{\sqrt{2}}\sqrt{1+3\cos^2\theta}, & n=3, \\
\sqrt{1+\cos^2\theta}, & n=4, \\
\tfrac{1}{2}(1-\sqrt{5})\cos\theta+\tfrac{1}{\sqrt{2}} \sqrt{3+\sqrt{5}+(5-\sqrt{5})\cos^2\theta}, & n=5;
\end{cases} \\
R_B
&=\begin{cases}
\tfrac{1}{\sqrt{2}}\sin\theta+\tfrac{1}{\sqrt{2}}\sqrt{1+3\sin^2\theta}, & n=3, \\
\sin\theta+\sqrt{2}\sqrt{1+\sin^2\theta}, & n=4, \\
\tfrac{1}{2}(1+\sqrt{5})\sin\theta+\tfrac{1}{\sqrt{2}}\sqrt{7+3\sqrt{5}+(5+\sqrt{5})\sin^2\theta}, & n=5;
\end{cases} \\
R_C
&=\begin{cases}
\dfrac{2\cos\theta\sin\theta-1}{\sqrt{2}(\cos\theta+\sin\theta)}, & n=3, \\
\dfrac{\cos\theta\sin\theta-\sqrt{2}}{\cos\theta+\sqrt{2}\sin\theta}, & n=4, \\
\dfrac{(\sqrt{5}-1)(\cos\theta\sin\theta-\sqrt{5}-2)}{2\cos\theta+(\sqrt{5}+1)\sin\theta}, & n=5.
\end{cases}
\end{align*}
The range for the angle $\theta$ is $[\frac{1}{2}\pi,\pi]$ for $\gamma_A$, $[\frac{3}{2}\pi,2\pi]$ for $\gamma_B$, and $[\pi,\frac{3}{2}\pi]$ for $\gamma_C$. 

Finally, we remark that for $\gamma_A$ and $\gamma_B$, $R$ is given by another quadratic equation
\[
R^2-2KR-L=0,\quad
R=\sqrt{K^2+L}+K,
\]
where $K$ and $L$ are given below. 

\begin{center}
\begin{tabular}{|c||c|c|c|}
\hline 
& $3$   
& $4$   
& $5$ \\
\hline  \hline 
$K_A$
& $\tfrac{1}{\sqrt{2}}\cos\theta$
& $0$
& $-\tfrac{1}{2}(\sqrt{5}-1)\cos\theta$ \\
\hline
$L_A$
& $\cos^2\theta+\tfrac{1}{2}$
& $\cos^2\theta+1$
& $\cos^2\theta+\tfrac{1}{2}(\sqrt{5}+3)$ \\
\hline 
$K_B$
& $\tfrac{1}{\sqrt{2}}\sin\theta$
& $\sin\theta$
& $\tfrac{1}{2}(\sqrt{5}+1)\sin\theta$ \\
\hline
$L_B$
& $\sin^2\theta+\tfrac{1}{2}$
& $\sin^2\theta+2$
& $\sin^2\theta+\tfrac{1}{2}(3\sqrt{5}+7)$ \\
\hline
\end{tabular}
\end{center}

\section{Features of Moduli Space} 
\label{feature}

\subsection{Companion}

The {\em companion} of pentagonal subdivision is introduced at the end of \cite[Section 3.1]{wy1}. In Figure \ref{companion}, the companion of the pentagon with anchor point $V$ is the pentagon with anchor point $V'$. We see that $V'$ is the reflection of $V$ with respect to the arc $BB'$ in Figure \ref{div_triangle}. For a pentagonal subdivision to have the companion, therefore, the necessary and sufficient condition is both $V$ and $V'$ are inside the moduli space. 

In the $M$-projection, we have the formulae for $R_B$ and $R_C$. We may verify that $R_B(2\pi-\rho)\ge R_C(\pi+\rho)$ for $\rho\in [0,\frac{1}{2}\pi]$. This means that the reflection $\gamma'_B$ of $\gamma_B$ with respect to the arc $BB'$ lies inside the moduli space. Therefore a pentagonal subdivision has companion, if and only if the anchor point $V$ lies in the region bounded by the arc $AA'$, the curve $\gamma_B$, and the curve $\gamma'_B$. 

In the $M$-projection, since the cartesian form of $\gamma_B$ is even in $x$, we know $\gamma'_B$ has the same algebraic formula as $\gamma_B$. In other words, the region for the companion to exist is bounded by the $x$-axis, and the $y\le 0$ part of the algebraic formula for $\gamma_B$.

\begin{figure}[htp]
\centering
\begin{tikzpicture}[>=latex,scale=1.5]

%% actual

\foreach \a in {0,1,2}
\foreach \b/\c in {0/1, 1/1, 2/-1}
{
\begin{scope}[xshift=3*\b cm, yscale=\c, rotate=120*\a]

\draw[gray]
	(0:1) arc (30:90:{sqrt(3)});

\coordinate (M) at (60:{sqrt(3)-1});
\coordinate (X) at (0:1);

\pgfmathsetmacro{\th}{18} 

\pgfmathsetmacro{\ph}{-atan((1/3)*(cot(\th+30)))}
\coordinate (A1) at ({\ph}:{sqrt(2)});
\coordinate (B1) at ({\ph+180}:{sqrt(2)});

\pgfmathsetmacro{\ps}{-atan((1/3)*(cot(\th+150)))}
\coordinate (A2) at ({\ps}:{sqrt(2)});
\coordinate (B2) at ({\ps+180}:{sqrt(2)});

\coordinate (V) at (36.8:0.595);
\coordinate (W) at (-42:0.95);
\coordinate (W2) at (120-42:0.95);

\arcThroughThreePoints[line width=1.2]{X}{A1}{V};
\arcThroughThreePoints[line width=1.2]{W}{A2}{X};
\arcThroughThreePoints[dashed]{V}{M}{W2};
	
\end{scope}
}

\foreach \a in {0,1,2}
{

\draw[rotate=120*\a] (0,0) -- (36.8:0.595);

\draw[xshift=3 cm, rotate=120*\a] (0,0) -- (-42:0.95);

\draw[xshift=6 cm, rotate=120*\a] (0,0) -- (42:0.95);

}

\foreach \b in {0,1,2}
{
\begin{scope}[xshift=3*\b cm]

\filldraw[fill=white]
	(0,0) circle (0.05);

\fill 
	(1,0) circle (0.05);

\end{scope}
}

\filldraw[fill=gray!50]
	(36.8:0.595) circle (0.05);
	
\node at (18:0.5) {$V$};

\filldraw[xshift=6 cm, yscale=-1, fill=gray!50]
	(-42:0.95) circle (0.05);

\node[xshift=9 cm] at (42:1.15) {$V'$};

\draw[->, very thick]
	(1.3,0) -- ++(0.5,0);

\draw[->, very thick]
	(4.3,0) -- ++(0.5,0);

\end{tikzpicture}
\caption{The companion of $V$ is $V'$.}
\label{companion}
\end{figure}

We remark that a pentagonal subdivision is the companion of itself, if and only if $V=V'$. This means exactly $V$ lies on $BM$.

\subsection{Location of Boundary}

In Figures \ref{projection-tetra}, \ref{projection-octa}, \ref{projection-icoA}, \ref{projection-icoB}, we divide the sphere into triangular regions. We know the locations of $\gamma_A,\gamma_B,\gamma_C$ in terms of these regions. However, the regions do not have the same size. The smallest regions are $\Omega_1,\Omega_2,\Omega_3,\Omega_4,\Omega_7,\Omega_8$, with area $(\frac{1}{n}-\frac{1}{6})\pi$. We further divide the sphere into triangular regions of the same size or even smaller, and then find more precise locations of $\gamma_A,\gamma_B,\gamma_C$ in terms of these smaller regions.

In the first of Figure \ref{boundary_curve}, we have the $A$-projection of one regular triangular face of the tetrahedron. We have the regions $\Omega_1,\Omega_2,\Omega_3,\Omega_4$ as before. Let $B',M'$ be the rotations of $B,M$ by $\frac{2}{3}\pi$. Let $B''$ be the rotation of $B$ by $\frac{4}{3}\pi$. Let $AX$ be the ray at angle $\frac{5}{6}\pi$, intersecting $B'B''$ at $X$. Then $\triangle AB'X$ is part of $\Omega_5$. 

We know $\gamma_A$ lies in $\Omega_5$, and between rays $AB'$ and $AX$. For any point \circled{$a$} on $B'X$, we get \squared{$a$} by rotating $-\frac{2}{3}\pi$ around $A$. Then \squared{$a$} lies on $BM$, and the arc $B$\circled{$a$} lies inside $\triangle BB'M'$. This implies that $A$\squared{$a$} intersects $B$\circled{$a$}. Therefore $B'X$ is outside the moduli space. Combined with the fact that $\gamma_A$ connects $A$ and $B'$, $AB'$ is in the moduli space, and $AX$ is not in the moduli space, we conclude that $\gamma_A$ is inside $\triangle AB'X$.

We also draw $\gamma_C$ in the first of Figure \ref{boundary_curve}. We know the $\Omega_4$-part of $\gamma_C$ is between the ray $AC$ and the ray $AY$ at angle $\frac{3}{2}\pi$. The point $Y$ is the intersection of $M'C$ and the ray $AY$. Then the picture shows that any point \circled{$c$} on $CY$ is not inside the moduli space. Combined with the fact that the $\Omega_4$-part of $\gamma_C$ connects $A$ and $C$, $AC$ is in the moduli space, and $AY$ is not in the moduli space, we conclude that the $\Omega_4$-part of $\gamma_C$ is inside $\triangle ACY$.

The $A$-projections for the octahedron and icosahedron are similar to the first of Figure \ref{boundary_curve}, and we have the similar $\triangle AB'X$ and $\triangle ACY$ for $n=4,5$. Since the argument above actually does not use the tetrahedron, we find $\gamma_A$ is inside $\triangle AB'X$, and the $\Omega_4$-part of $\gamma_C$ is inside $\triangle ACY$, also for $n=4,5$.

\begin{figure}[htp]
\centering
\begin{tikzpicture}[>=latex]

% A-projection n=3

\begin{scope}[scale=3]

\pgfmathsetmacro{\dabt}{sqrt(3)-1};
\pgfmathsetmacro{\dabtt}{2/(sqrt(3)-1)};

\coordinate (M1) at ({-\dabt},0);
\coordinate (MM1) at ({\dabtt},0);
\coordinate (C1) at (-60:{\dabt});
\coordinate (CC1) at (120:{\dabtt});
\coordinate (M) at (60:{\dabt});

\coordinate (A1) at (1,0);
\coordinate (X1) at (140:0.84);
\coordinate (XX1) at (20:0.84);
\coordinate (AA1) at (-2,0);

\coordinate (Y1) at (165:0.59);
\coordinate (YY1) at (285:0.59);
 
\foreach \a in {0,1,2}
\draw[rotate=120*\a,gray]
	(1,0) arc (30:90:{sqrt(3)})
	(120:1) -- (-60:{\dabt});
	
\draw[gray]
	(0,0) -- (150:0.79)
	(0,0) -- (270:0.79);

\fill[gray!50]
	(150:0.79) -- (0,0) -- (120:1) arc (150:170:{sqrt(3)})
	(0,-0.504) -- (0,0) -- (-60:{\dabt}) -- (285:0.59);

\arcThroughThreePoints[gray]{M1}{MM1}{C1};
\arcThroughThreePoints[gray]{M}{MM1}{M1};
%\arcThroughThreePoints[red]{M1}{C1}{MM1};
\arcThroughThreePoints{A1}{AA1}{X1};\arcThroughThreePoints{Y1}{CC1}{C1};

\draw
	(0,0) -- (X1)
	(0,0) -- (XX1)
	(0,0) -- (Y1)
	(0,0) -- (YY1);

\draw[densely dotted, scale=sqrt(2)]
	plot [domain=60:89.5, samples=100] ({60+\x}:{sqrt(1+(1/32)*sec(\x)*sec(\x))-(1/(4*sqrt(2)))*sec(\x)})
	plot [domain=-89.5:0, samples=100] ({\x}:{sqrt(1+(1/8)*sec(\x)*sec(\x))-(1/(2*sqrt(2)))*sec(\x)});
		
\fill 
	(1,0) circle (0.033);
\filldraw[fill=white]
	(0,0) circle (0.033);
	
\node[draw, fill=white, shape=circle, inner sep=0.3] at (X1) {\tiny $a$};
\node[draw, fill=white, shape=rectangle, inner sep=1] at (XX1) {\tiny $a$};
\node[draw, fill=white, shape=circle, inner sep=0.3] at (YY1) {\tiny $c$};
\node[draw, fill=white, shape=rectangle, inner sep=1] at (Y1) {\tiny $c$};
	
\node at (0.1,-0.05) {\small $A$};
\node at (1.1,0) {\small $B$};
\node at (120:1.05) {\small $B'$};
\node at (240:1.1) {\small $B''$};
\node at (-60:0.8) {\small $C$};
\node at (60:0.8) {\small $M$};
\node at (-0.83,0) {\small $M'$};
\node at (150:0.87) {\small $X$};
\node at (-0.05,-0.55) {\small $Y$};

\node at (-0.45,0.65) {\small $\gamma_A$};
\node at (0.87,-0.55) {\small $\gamma_C$};

\node at (0.5,0.4) {\small $1$};
\node at (0.5,-0.4) {\small $2$};
\node at (-0.2,0.7) {\small $3$};
\node at (-0.2,-0.7) {\small $4$};
\node at (-0.85,0.6) {\small $5$};
\node at (-0.85,-0.6) {\small $6$};

\end{scope}

% B-projection, n=4

\begin{scope}[xshift=9.5cm, scale=3.7]

\coordinate (A2) at ({-sqrt(3)+1},0);
\coordinate (AA2) at ({sqrt(3)+1},0);
\coordinate (C2) at (45:{sqrt(2)-2});
\coordinate (C2A) at (45:{2+sqrt(2)});
\coordinate (CC2) at (-45:{2-sqrt(2)});
\coordinate (CCC2) at (45:{2-sqrt(2)});
\coordinate (MM2) at ({-sqrt(2)},0);
\coordinate (MMM2) at (0,{-sqrt(2)});

\coordinate (X2) at (65:0.61);
\coordinate (XX2) at (155:0.61);

\coordinate (Y2) at (-25:0.486);
\coordinate (YY2) at (-115:0.486);

\foreach \b in {0,1,2,3}
\draw[gray,rotate=90*\b]
	(45:{2-sqrt(2)}) -- (0,0) -- ({sqrt(3)-1},0) arc (30:60:2);

\draw[gray]
	(MM2) -- (A2);
	
\arcThroughThreePoints[gray]{MM2}{C2}{CC2};
\arcThroughThreePoints[gray]{CC2}{MMM2}{CCC2};

\fill[gray!50]
	(0,{sqrt(3)-1}) -- (0,0) -- (45:{2-sqrt(2)}) to[out=135,in=-30] (0,{sqrt(3)-1})
	(0,-0.45) -- (0,0) -- (45:{sqrt(2)-2})
to[out=-9,in=180] (0,-0.45);
	
\draw
	(0,0) -- (X2)
	(0,0) -- (XX2)
	(0,0) -- (Y2)
	(0,0) -- (YY2);
	
\arcThroughThreePoints{X2}{AA2}{A2};
\arcThroughThreePoints{C2}{C2A}{Y2};	
	
\draw[densely dotted, scale=sqrt(2)]
	plot [domain=90.5:135, samples=100] ({-45+\x}:{sqrt(1+(1/4)*sec(\x)*sec(\x))+(1/2)*sec(\x)})
	plot [domain=-90.5:-180, samples=100] ({\x}:{sqrt(1+(1/2)*sec(\x)*sec(\x))+(1/sqrt(2))*sec(\x)})
	;

\fill
	(0,0) circle (0.025);
	
\filldraw[fill=white]
	({-sqrt(3)+1},0) circle (0.025);

\node[draw, fill=white, shape=circle, inner sep=0.3] at (X2) {\tiny $b$};
\node[draw, fill=white, shape=rectangle, inner sep=1] at (XX2) {\tiny $b$};
\node[draw, fill=white, shape=circle, inner sep=0.3] at (YY2) {\tiny $c$};
\node[draw, fill=white, shape=rectangle, inner sep=1] at (Y2) {\tiny $c$};

\node at (-0.8,0) {\small $A$};
\node at (0,0.78) {\small $A_1$};
\node at (0.08,0.02) {\small $B$};
\node at (45:-0.65) {\small $C$};
\node at (-45:0.65) {\small $C_1$};
\node at (135:0.65) {\small $M$};
\node at (45:0.65) {\small $M_1$};
\node at (-1.48,0) {\small $M'$};

\node at (0.07,0.57) {\small $\gamma_B$};
\node at (-0.75,-0.22) {\small $\gamma_C$};

\node at (-0.43,0.3) {\small $1$};
\node at (-0.43,-0.3) {\small $2$};
\node at (-0.1,0.5) {\small $7$};
\node at (-0.1,-0.5) {\small $8$};
\node at (0.3,0.7) {\small $13$};
\node at (0.3,-0.7) {\small $14$};

\end{scope}

% B-projection, n=5

\begin{scope}[shift={(5.5cm,-6cm)}, scale=12]

\pgfmathsetmacro{\ra}{sqrt((sqrt(5)+5)/2)};
\pgfmathsetmacro{\dab}{(sqrt(30+6*sqrt(5))-sqrt(5)-3)/4};
\pgfmathsetmacro{\dbm}{(sqrt(10+2*sqrt(5))-sqrt(5)-1)/2};

\coordinate (A5) at (-\dab,0);
\coordinate (A51) at (108:\dab);
\coordinate (A52) at (36:\dab);
\coordinate (A53) at (-36:\dab);
\coordinate (AA5) at ({1/\dab},0);
\coordinate (AA51) at (-72:{1/\dab});	

\coordinate (C5) at (216:\dbm);
\coordinate (C51) at (-72:\dbm);
\coordinate (CC5) at (36:{1/\dbm});
\coordinate (CC51) at (108:{1/\dbm});
\coordinate (M5) at (144:\dbm);
\coordinate (M51) at (72:\dbm);
\coordinate (M52) at (0:\dbm);

\coordinate (MP) at (-0.557536516,0);
\coordinate (MMP) at (1.793604493,0);	
\coordinate (M1P) at (72:-0.557536516);
\coordinate (MM1P) at (72:1.793604493);

\coordinate (B5) at (85:0.291);
\coordinate (BB5) at (157:0.291);
\coordinate (D5) at (62:0.1297);
\coordinate (DD5) at (134:0.1297);
\coordinate (E5) at (0.015,0.1735);	

\coordinate (U5) at (240:0.2407);
\coordinate (UU5) at (-48:0.2407);
\coordinate (V5) at (240:0.14);
\coordinate (VV5) at (-48:0.14);
\coordinate (W5) at (260:0.107);
\coordinate (WW5) at (-28:0.107);

\arcThroughThreePoints[gray]{A52}{AA5}{A5};	
\arcThroughThreePoints[gray]{A53}{AA51}{A51};
\arcThroughThreePoints[gray]{C5}{CC5}{M52};	
\arcThroughThreePoints[gray]{C51}{CC51}{M51};	
\arcThroughThreePoints[gray]{MP}{MMP}{C51};	
\arcThroughThreePoints[gray]{C51}{M1P}{M52};

\foreach \b in {0,...,4}
\draw[rotate=72*\b,gray]
	(144:{\dbm}) -- (0,0) -- (-\dab,0) arc (150:138:\ra);

\fill[gray!50]
	(A51) -- (100:0.252) -- (85:0.175) -- (54:0.122) -- (0,0) -- (M51) to[out=162,in=-12] (A51)
	(0,0) -- (0,-0.0995) -- (V5) -- (C5) -- (U5) -- (-108:0.2354);

\draw[gray]
	(0,0) -- (MP)
	(0,0) -- (54:0.297)
	(0,0) -- (-90:0.297);
	
\draw
	(0,0) -- (B5)
	(0,0) -- (BB5)
	(0,0) -- (D5)
	(0,0) -- (DD5)
	(0,0) -- (U5)
	(0,0) -- (UU5)
	(0,0) -- (W5)
	(0,0) -- (WW5)
	;
	
\arcThroughThreePoints{B5}{AA5}{A5};	
\arcThroughThreePoints{D5}{AA5}{A5};
\arcThroughThreePoints{E5}{AA5}{A5};
\arcThroughThreePoints{C5}{CC5}{UU5};
\arcThroughThreePoints{C5}{CC5}{VV5};
\arcThroughThreePoints{C5}{CC5}{WW5};

\draw[densely dotted]
	plot [domain=90.5:144, samples=100] ({-36+\x}:{sqrt(1+1.1215169944*sec(\x)*sec(\x))+1.0590169944*sec(\x)})
	plot [domain=-91:-180, samples=100] ({\x}:{sqrt(1+1.7135*sec(\x)*sec(\x))+1.3090*sec(\x)})
	;

\node[draw, fill=white, shape=circle, inner sep=0.3] at (B5) {\tiny $a$};
\node[draw, fill=white, shape=rectangle, inner sep=1] at (BB5) {\tiny $a$};
\node[draw, fill=white, shape=circle, inner sep=0.3] at (85:0.175) {\tiny $b$};
\node[draw, fill=white, shape=rectangle, inner sep=1] at (157:0.175) {\tiny $b$};
\node[draw, fill=white, shape=circle, inner sep=0.3] at (D5) {\tiny $c$};
\node[draw, fill=white, shape=rectangle, inner sep=1] at (DD5) {\tiny $c$};

\node[draw, fill=white, shape=circle, inner sep=0.3] at (U5) {\tiny $u$};
\node[draw, fill=white, shape=rectangle, inner sep=1] at (UU5) {\tiny $u$};
\node[draw, fill=white, shape=circle, inner sep=0.3] at (V5) {\tiny $v$};
\node[draw, fill=white, shape=rectangle, inner sep=1] at (VV5) {\tiny $v$};
\node[draw, fill=white, shape=circle, inner sep=0.3] at (W5) {\tiny $w$};
\node[draw, fill=white, shape=rectangle, inner sep=1] at (WW5) {\tiny $w$};

\fill
	(0,0) circle (0.009);
	
\filldraw[fill=white]
	(-\dab,0) circle (0.009);

\node at (-0.36,0) {\small $A$};
\node at (108:0.36) {\small $A_1$};
\node at (36:0.36) {\small $A_2$};
\node at (-36:0.36) {\small $A_3$};
\node at (-108:0.36) {\small $A_4$};
\node at (0.025,0.01) {\small $B$};
\node at (38:-0.3) {\small $C$};
\node at (144:0.3) {\small $M$};
\node at (72:0.3) {\small $M_1$};
\node at (-72:0.305) {\small $C_1$};
\node at (0:0.31) {\small $M_2$};
\node at (-0.58,0) {\small $M'$};
\node at (47:0.13) {\small $X$};
\node at (66:0.16) {\small $N$};
\node at (0.017,-0.1) {\small $Y$};

\node at (0.05,0.21) {\small $\gamma_B$};
\node at (-0.35,-0.07) {\small $\gamma_C$};

\node at (-0.14,0.03) {\small 1};
\node at (-0.14,-0.03) {\small 2};
\node at (-0.12,0.25) {\small 7};
\node at (-0.12,-0.25) {\small 8};
\node at (0,0.33) {\small 13};
\node at (0,-0.33) {\small 14};
\node at (0.2,0.27) {\small 19};
\node at (0.2,-0.27) {\small 20};

\end{scope}

\end{tikzpicture}
\caption{Locations of the boundary curves.}
\label{boundary_curve}
\end{figure}

For the location of $\gamma_B$ and the remaining part of $\gamma_C$, we use the $B$-projection. For the tetrahedron ($n=3$), the picture is the the same as the horizontal flip of the $A$-projection (the first of Figure \ref{boundary_curve}), with $A$ and $B$ exchanged. In fact, $\gamma_A$ and the $\Omega_4$-part of $\gamma_C$ are also transformed to $\gamma_B$ and the $\Omega_8$-part of $\gamma_C$ in this way. Then we get the corresponding triangles containing $\gamma_B$ and the $\Omega_8$-part of $\gamma_C$.  

The second of Figure \ref{boundary_curve} is the $B$-projection for the octahedron. Let $A_1,M_1$ be the rotations of $A,M$ around $B$ by $-\frac{1}{2}\pi$ ($A_1$ is $A'$ in Figure \ref{3curves}). Let $C_1$ be the rotation of $C$ by $\frac{1}{2}\pi$. Then $A_1M_1$ cuts a triangle from $\Omega_{19}$, and $CC_1$ cuts a triangle from $\Omega_8$. Both triangles are indicated by gray shade. 

We know $\gamma_B$ connects $B$ and $A_1$, and $BM_1$ is not in the moduli space. The picture shows that any point \circled{$b$} on $A_1M_1$ is not inside the moduli space. Therefore $\gamma_B$ is inside $\triangle A_1BM_1$. Similar argument shows that the $\Omega_8$-part of $\gamma_C$ lies in the other shaded triangle. We also note that $M',C,C_1$ lies in the same great circle. Therefore the whole $M'C_1$ is outside $\gamma_C$. Since $\gamma_C$ passes $C$, we conclude that $M'C_1$ is tangential to $\gamma_C$.

The third of Figure \ref{boundary_curve} is the $B$-projection for the icosahedron. Let $A_k,M_k$ be the rotations of $A,M$ by $-\frac{2k}{5}\pi$. Let $C_1$ be the rotation of $C$ by $\frac{2}{5}\pi$. Let $BX$ and $BY$ be the rays at angles $\frac{3}{10}\pi$ and $\frac{3}{2}\pi$. We see $BX$ intersects $A_1A_3$ at $X$, and $BY$ intersects $CM_2$ at $Y$. 

The arcs $AA_2,A_1A_3,BM_1$ intersect at a point $N$. It is easy to verify that $N$ is a point on $\gamma_B$. By Lemma \ref{geometry2}, this implies $BN$ is inside the moduli space, and $NM_1$ is outside the moduli space. We also know $BX$ is outside the moduli space. The picture shows that $A_1N$ is inside the moduli space, and $A_1M_1,NX$ are outside the moduli space. Then we conclude $\gamma_B$ is inside $\triangle A_1M_1N\cup\triangle BNX$.

The arcs $CM_2,BA_4$ intersect at a point, that we can easily verify to be on $\gamma_C$. We have two shaded triangles bounded by $BA_4,BY,CC_1,CM_2$. By an argument similar to $\gamma_B$, we know $\gamma_C$ is inside the union of the two shaded triangles. We also note that, similar to the octahedron, $M',C,C_1$ are on the same great circle. This implies $M'C_1$ is tangential to $\gamma_C$ at $C$.

\subsection{Upper Bound of Edge Length}

\begin{theorem}
The suprema of $a,b,c$ in the pentagonal subdivision tiling are given by Table \ref{supremum}. 
\end{theorem}

\begin{table}[htp]
\centering
\begin{tabular}{|c||c|c|c|}
\hline 
& $n=3$   
& $n=4$  
& $n=5$  \\
\hline  \hline 
$\sup a$
& $0.6259072384\pi$
& $0.4171042591\pi$
& $0.2569674158\pi$ \\
\hline
$\sup b$
& $0.6259072384\pi$
& $\tfrac{1}{2}\pi$
& $\arccos\tfrac{1}{\sqrt{5}}$ \\
\hline
$\sup c$
& $\pi$
& $0.6882086654\pi$
& $0.4357713623\pi$ \\
\hline 
\end{tabular}
\caption{The suprema of the edges in the pentagonal subdivision tiling.}
\label{supremum}
\end{table}

The decimal value $0.6259072384\pi$ means 
\[
0.6259072383\pi<\sup a< 0.6259072384\pi.
\]
In fact, we obtain the precise value 
\[
\sup a=2\arctan (\sqrt{R_B(\theta_{AB})^2+1}-R_B(\theta_{AB})),
\]
where $R_B(\theta)$ is given by \eqref{rb}, and the precise value of $\tan\theta_{AB}$ is given after \eqref{rb}, in terms of square and cube roots.

The same remark applies to the other decimal values, and we also obtain the precise values. 

The supremum of $a$ is the biggest of the maximal distances from $A$ to $\gamma_A,\gamma_B,\gamma_C$. Here we have the supremum instead of the maximum because the moduli space is an open subset.

In the $A$-projection, the bigger Euclidean distance from $A$ (this is $r$ in \eqref{eqE} and \eqref{eqD}) is equivalent to the bigger spherical distance from $A$. We also know the polar versions \eqref{eqE} and \eqref{eqD} of $\gamma_A$ and $\gamma_C$ are strictly monotone with respect to the angle $\phi$. Moreover, the whole $\gamma_C$ is given by the same formulae \eqref{eqE} and \eqref{eqD}, and the whole range is $\phi\in [\alpha,\alpha+\frac{1}{2}\pi]$. Then we find $R$ and $r$ are strictly monotone in $\phi$ on the whole range. Therefore we get (see Figure \ref{3curves})
\[
d_A(A,\gamma_A)=|AB'|_A=d_{AB},\quad
d_A(A,\gamma_C)=|AB|_A=d_{AB}.
\]
Here $d_A(A,\gamma)$ is the Euclidean distance from $A$ to curve $\gamma$ in the $A$-projection. By applying the same reason to the $B$-projection, we get
\[
d_B(B,\gamma_B)=|BA'|_B=d_{AB},\quad
d_B(B,\gamma_C)=|BA|_B=d_{AB}.
\]
In Figure \ref{boundary_curve}, we always have $d_A(A,\gamma_B)\ge |AA_1|_A=|AA'|_A>d_{AB}$. We also have $d_B(B,\gamma_A)\ge |BB'|_B>d_{AB}$. Therefore 
\[
\sup a=d_{\bb S}(A,\gamma_B),\quad
\sup b=d_{\bb S}(B,\gamma_A).
\]
Here the subscript ${\bb S}$ indicates the spherical distance.

For $n=4$, we have (the spherical distance) $|BP|_{\bb S}=|BB'|_{\bb S}=\frac{1}{2}\pi$ for any $P$ on $B'M'$. For $n=5$, we know $|BP|_{\bb S}$ is strictly decreasing as $P$ moves from $B'$ to $M'$ along $B'M'$. Therefore
\[
\sup b=|BB'|_{\bb S}
=\begin{cases}
\frac{1}{2}\pi, & n=4, \\
\arccos\frac{1}{\sqrt{5}}, & n=5.
\end{cases}
\]
By the symmetry, we have $\sup a=\sup b$ for $n=3$. Therefore for the bounds of $a$ and $b$, it remains to calculate $d_{\bb S}(A,\gamma_B)$.

As to the supremum of $c$, we note that $\sup \frac{1}{2}c$ is the biggest of the maximal distances from $M$ to $\gamma_A,\gamma_B,\gamma_C$. The three maximal distances can be calculated by the $M$-projections of the three curves. This corresponds to minimizing $R_A,R_B,R_C$ given at the end of Section \ref{picture}. 

For $\min R_C$, we recall the range $\theta\in [\pi,\frac{3}{2}\pi]$ for $\gamma_C$. For $n=3$, the solution $\theta_{MC}$ of $R_C'(\theta)=0$ satisfies
\[
\tan\theta_{MC}=1.
\]
Correspondingly, we have 
\[
\theta_{MC}=\tfrac{5}{4}\pi,\;
R_C(\theta_{MC})=0.
\]
By $R_C(\theta_{MC})<R_C(\pi)=R_C(\frac{3}{2}\pi)=\tfrac{1}{\sqrt{2}}$, we get 
\[
\min R_C
=R_C(\theta_{MC}).
\]
Correspondingly, we have 
\begin{align*}
\max r_C
&=r_C(\theta_{MC})
={\textstyle \sqrt{R_C(\theta_{MC})^2+1}}-R_C(\theta_{MC})=1, \\
d_{\bb S}(M,\gamma_C)
&=2\arctan \max r_C=\tfrac{1}{2}\pi.
\end{align*}

For $n=4$, the solution $\theta_{MC}$ of $R_C'(\theta)=0$ satisfies
\[
\tan\theta_{MC}=\tfrac{1}{6\sqrt{2}}((\sqrt{17}-5)T^{\frac{1}{3}}+(\sqrt{17}+5)T^{-\frac{1}{3}}+2),\quad
T=\tfrac{1}{4}(\sqrt{17}+1).
\]
The corresponding 
\begin{align*}
R_C(\theta_{MC})
&=\tfrac{\sqrt{17}}{18\sqrt{3}}\left( 
(\sqrt{17}-1)T^{\frac{1}{3}}
-(\sqrt{17}+1)T^{-\frac{1}{3}}
+\tfrac{10}{\sqrt{17}} \right) \\
&\quad \sqrt{ 
(\sqrt{17}+1)T^{\frac{1}{3}}
+(\sqrt{17}-1)T^{-\frac{1}{3}}+5 } \\
&=\tfrac{\sqrt{17}}{9\sqrt{3}}\left( 
2T^{-\frac{2}{3}}
-2T^{\frac{2}{3}}
+\tfrac{5}{\sqrt{17}} \right)
\sqrt{ 
4T^{\frac{4}{3}}
+4T^{-\frac{4}{3}}+5 } \\
&=0.53308.
\end{align*}
We may verify $R_C(\theta_{MC})<R_C(\pi)=R_C(\frac{3}{2}\pi)$. This implies $\min R_C=R_C(\theta_{MC})$.

For $n=5$, the solution $\theta_{MC}$ of $R_C'(\theta)=0$ satisfies
\begin{align*}
\tan\theta_{MC}
&=\tfrac{1}{6\sqrt{5}\sqrt{1+3\sqrt{5}}}
\left((-\sqrt{3}(11+7\sqrt{5}){\textstyle \sqrt{\sqrt{5}-2}}+\sqrt{14}{\textstyle \sqrt{7+9\sqrt{5}}})T^{\frac{1}{3}}\right. \\
&\qquad\qquad\quad 
+(\sqrt{3}(11+7\sqrt{5}){\textstyle \sqrt{\sqrt{5}-2}}+\sqrt{14}{\textstyle \sqrt{7+9\sqrt{5}}})T^{-\frac{1}{3}} \\
&\qquad\qquad\quad \left.
+(3+\sqrt{5}){\textstyle \sqrt{1+3\sqrt{5}}} \right), \\
T
&=\frac{ 3\sqrt{3}(3-\sqrt{5})\sqrt{2+\sqrt{5}}+ \sqrt{70}\sqrt{7+9\sqrt{5}} }{2(1+3\sqrt{5})^{\frac{3}{2}}}.
\end{align*}
The corresponding 
\begin{align*}
R_C(\theta_{MC})
&=\tfrac{1}{60\sqrt{30}\sqrt{4+\sqrt{5}}(1+3\sqrt{5})^{\frac{1}{4}}} \\
&\quad \left((\sqrt{3}(17+21\sqrt{5}){\textstyle \sqrt{3+\sqrt{5}}}-\sqrt{140}{\textstyle \sqrt{19-\sqrt{5}}} )T^{\frac{1}{3}}\right. \\
&\quad
+(-\sqrt{3}(17+21\sqrt{5}){\textstyle \sqrt{3+\sqrt{5}}}-\sqrt{140}{\textstyle \sqrt{19-\sqrt{5}}})T^{-\frac{1}{3}} \\
&\quad \left.
+2(37-\sqrt{5}){\textstyle \sqrt{4+\sqrt{5}}}\right) \\
&\quad \left( 
 4( \sqrt{3}(\sqrt{5}-1){\textstyle \sqrt{\sqrt{5}-2}} + \sqrt{14}{\textstyle \sqrt{7+9\sqrt{5}}} )T^{\frac{1}{3}}\right. \\
&\quad +4(-\sqrt{3}(\sqrt{5}-1){\textstyle \sqrt{\sqrt{5}-2}} + \sqrt{14}{\textstyle \sqrt{7+9\sqrt{5}}} )T^{-\frac{1}{3}} \\
&\quad \left. +(27+\sqrt{5}){\textstyle \sqrt{1+3\sqrt{5}}} \right)^{\frac{1}{2}} \\
&=1.22527.
\end{align*}
We may also verify $R_C(\theta_{MC})<R_C(\pi)=R_C(\frac{3}{2}\pi)$. This implies $\min R_C=R_C(\theta_{MC})$.

Similarly, by the formulae for $R_A,R_B$, we may easily find
\[
\min R_A
=\min R_B
=\begin{cases}
\tfrac{1}{\sqrt{3}}=0.57735, & n=3, \\
1, & n=4, \\
\tfrac{1}{2}(\sqrt{5}+1)=1.61803, & n=5.
\end{cases}
\]
We always have $\min R_A=\min R_B>\min R_C$. This means 
\[
\max r_A
=\max r_B
<\max r_C
=r_C(\theta_{MC})
={\textstyle \sqrt{R_C(\theta_{MC})^2+1}}-R_C(\theta_{MC}).
\]
Then we have
\[
\sup c=2d_{\bb S}(M,\gamma_C)
=4\arctan r_C(\theta_{MC})
=\begin{cases}
\pi, & n=3, \\
0.6882086653\pi, & n=4, \\
0.4357713622\pi, & n=5.
\end{cases}
\]
The value $0.6882086653\pi$ means $0.6882086653\pi<\sup c< 0.6882086654\pi$. We add $10^{-10}\pi$ to get the values in Table \ref{supremum}.

Now we calculate $\sup a=d_{\bb S}(A,\gamma_B)$ in the similar way. First, we use the M\"obius transform to find the polar formula for the $A$-projection of $\gamma_B$. It is given by $r_B=\sqrt{R_B^2+1}-R_B$ with
\begin{equation}\label{rb}
R_B =\begin{cases}
\dfrac{\sqrt{3}-8\cos\theta\sin\theta}{2\sqrt{2}(\sqrt{3}\cos\theta+\sin\theta)}, & n=3, \\
\dfrac{\sqrt{3}-2\cos\theta\sin\theta}{\sqrt{2}(\sqrt{3}\cos\theta+\sin\theta)}, & n=4, \\
\dfrac{\sqrt{3}(3\sqrt{5}+7)-8\cos\theta\sin\theta}{2(3-\sqrt{5})(\sqrt{3}\cos\theta+\sin\theta)}, & n=5.
\end{cases}
\end{equation}
The range is $\theta\in [0,\frac{1}{3}\pi]$. For $n=3$, the solution $\theta_{AB}$ of $R_B'(\theta)=0$ satisfies
\begin{align*}
\tan\theta_{AB}
&=\tfrac{1}{77\sqrt{3}}((15\sqrt{3}+\sqrt{451})T^{\frac{1}{3}}+7(-15\sqrt{3}+\sqrt{451})T^{-\frac{1}{3}}+7), \\
T &=6\sqrt{3}+\sqrt{451}.
\end{align*}
For $n=4$, we have
\begin{align*}
\tan\theta_{AB}
&=\tfrac{1}{25\sqrt{3}}((6\sqrt{3}+\sqrt{38})T^{\frac{1}{3}}+5(-6\sqrt{3}+\sqrt{38})T^{-\frac{1}{3}}+5), \\
T &=3\sqrt{3}+2\sqrt{38}.
\end{align*}
For $n=5$, we have
\begin{align*}
\tan\theta_{AB}
&=\tfrac{1}{218\sqrt{6}(7+4\sqrt{5})}
\left((\sqrt{2}(5\sqrt{5}-4){\textstyle \sqrt{501+127\sqrt{5}}}
+3\sqrt{3}(31+43\sqrt{5}))T^{\frac{1}{3}}\right. \\
&\left. 
+(27+11\sqrt{5})(\sqrt{2}(5\sqrt{5}-4){\textstyle \sqrt{501+127\sqrt{5}}}
-3\sqrt{3}(31+43\sqrt{5}))T^{-\frac{1}{3}} \right)  \\
&+\tfrac{1}{109\sqrt{3}}(17+6\sqrt{5}), \\
T &=(3\sqrt{5}+7)(3\sqrt{6}+{\textstyle \sqrt{501+127\sqrt{5}}}
).
\end{align*}
We verify the corresponding $\theta_{AB}$ is always in $[0,\frac{1}{3}\pi]$. Then we verify $R_B(\theta_{AB})<R_B(\frac{1}{3}\pi)<R_B(0)$ always holds. Then we conclude
\[
\sup a
=d_{\bb S}(A,\gamma_B)
=2\arctan r_B(\theta_{AB})
=\begin{cases}
0.6259072383\pi, & n=3, \\
0.4171042590\pi, & n=4, \\
0.2569674157\pi, & n=5.
\end{cases}
\]
We add $10^{-10}\pi$ to get the values in Table \ref{supremum}.

\subsection{Reduction}

The pentagon in the pentagonal subdivision of a Platonic solid has three edge lengths $a,b,c$. The edge lengths are generally not equal, and we get a tiling of the sphere by congruent pentagons with the edge length combination $a^2b^2c$. When some edge lengths become equal, then the tiling has a different edge length combination. 

For example, if $a=b$, then the edge length combination $a^2b^2c$ is reduced to $a^4b$ (with $a=b$ labeled as the new $a$, and $c$ labeled as the new $b$), and we get a tiling of the sphere by congruent {\em almost equilateral} pentagons. The reduction means $|AV|_{\bb S}=|BV|_{\bb S}$. All the points of equal distance from $A$ and $B$ form a great circle. By $\tan\frac{1}{2}|AV|_{\bb S}=|z_A|$ and $\tan\frac{1}{2}|BV|_{\bb S}=|z_B|$, the equality $|AV|_{\bb S}=|BV|_{\bb S}$ becomes $|z_A|=|z_B|$. Let $z=z_M$ be the complex coordinate of $V$ in the $M$-projection. Then we substitute \eqref{m2a} and \eqref{m2b} into $|z_A|=|z_B|$, and get the equation of the circle
\begin{align*}
n=3 &\colon x-y=0, \\
n=4 &\colon x^2+y^2+2(\sqrt{3}+2)(\sqrt{2}x-\sqrt{3}y)-1=0, \\
n=5 &\colon 
x^2+y^2
-2((\sqrt{5}-11)\mu-\sqrt{5}) x 
-2((5\sqrt{5}+3)\mu+3)y
-1=0, \\
\mu &=\tfrac{1}{58}{\textstyle \sqrt{195-6\sqrt{5}}  }.
\end{align*}
For $n=3$, the circle is actually the diagonal line. For $n=4$, the circle has center $(-\sqrt{2}(2+\sqrt{3}),\sqrt{3}(2+\sqrt{3}))$ and radius $2\sqrt{9+5\sqrt{3}}$. For $n=5$, the center is $((\sqrt{5}-11)\mu-\sqrt{5}, (5\sqrt{5}+3)\mu+3)$, and the radius is $\sqrt{4(13\sqrt{5}+2)\mu+30}$.

The circle is also the intersection of the sphere with a plane. The following is the equation of this plane
\begin{align*}
n=3 &\colon \xi_1-\xi_2=0, \\
n=4 &\colon \sqrt{2}\xi_1-\sqrt{3}\xi_2-(\sqrt{3}-2)\xi_3=0, \\
n=5 &\colon 
2\cdot 5^{\frac{1}{4}}\xi_1-\sqrt{6}{\textstyle \sqrt{\sqrt{5}+1}}\xi_2+((\sqrt{5}+3)5^{\frac{1}{4}}-\tfrac{\sqrt{3}}{\sqrt{2}}(\sqrt{5}+1)^{\frac{3}{4}})\xi_3=0.
\end{align*}

If $a=c$ or $b=c$, then we get a tiling of the sphere by congruent pentagons with the edge length combination $a^3b^2$. The reduction $a=c$ is the same as $a=2c_1$, or $|AV|_{\bb S}=2|MV|_{\bb S}$. By $\tan\frac{1}{2}|AV|_{\bb S}=|z_A|$ and $\tan\frac{1}{2}|MV|_{\bb S}=|z_M|=|z|$, the equality $|AV|_{\bb S}=2|MV|_{\bb S}$ becomes 
\[
|z_A|=\frac{2|z|}{1-|z|^2}.
\]
Applying \eqref{m2a}, we get the equation in the $M$-projection
\begin{align*}
n=3 &\colon (x^2+y^2)^2
+\sqrt{2}(\sqrt{3}-1)(x^2+y^2)x
-3\sqrt{3}(\sqrt{3}-1)(x^2+y^2) \\
&\qquad 
+\sqrt{2}(\sqrt{3}-1)x
-\sqrt{3}+2=0, \\
n=4 &\colon (x^2+y^2)^2
+2(\sqrt{3}-\sqrt{2})(x^2+y^2)x
-6\sqrt{3}(\sqrt{3}-\sqrt{2})(x^2+y^2) \\
&\qquad 
+2(\sqrt{3}-\sqrt{2})x
-2\sqrt{6}+5=0, \\
n=5 &\colon 
(x^2+y^2)^2
+(\sqrt{5}-\sqrt{3})(\sqrt{3}-1)(x^2+y^2)x \\
&\qquad 
+3(2\sqrt{15}-3\sqrt{5}+4\sqrt{3}-9)(x^2+y^2)  \\
&\qquad
+(\sqrt{5}-\sqrt{3})(\sqrt{3}-1)x 
-2\sqrt{15}+3\sqrt{5}-4\sqrt{3}+8=0.
\end{align*}
It is very easy to convert the cartesian form to the polar form
\begin{align*}
n=3 &\colon r^4-3\sqrt{3}(\sqrt{3}-1)r^2-\sqrt{3}+2+\sqrt{2}(\sqrt{3}-1)r(r^2+1)\cos\theta=0 \\
n=4 &\colon r^4
-6\sqrt{3}(\sqrt{3}-\sqrt{2})r^2-2\sqrt{6}+5
+2(\sqrt{3}-\sqrt{2})r(r^2+1)\cos\theta
=0, \\
n=5 &\colon 
r^4+3(2\sqrt{15}-3\sqrt{5}+4\sqrt{3}-9)r^2
-2\sqrt{15}+3\sqrt{5}-4\sqrt{3}+8   \\
&\qquad
+(\sqrt{5}-\sqrt{3})(\sqrt{3}-1)r(r^2+1)\cos\theta=0.
\end{align*}
The curve is also the intersection of the sphere with a parabolic cylinder
\begin{align*}
n=3 &\colon 2\sqrt{3}\xi_3^2+\sqrt{2}\xi_1+\xi_3-\sqrt{3}=0, \\
n=4 &\colon 2\sqrt{3}\xi_3^2+\xi_1+\sqrt{2}\xi_3-\sqrt{3}=0, \\
n=5 &\colon 
4\sqrt{3}\xi_3^2+(\sqrt{5}-1)\xi_1+(\sqrt{5}+1)\xi_3-2\sqrt{3}=0.
\end{align*}

Similarly, the reduction $b=c$ means $|BV|_{\bb S}=2|MV|_{\bb S}$, or
\[
|z_B|=\frac{2|z|}{1-|z|^2}.
\]
Applying \eqref{m2b}, we get the equation in the $M$-projection
\begin{align*}
n=3 &\colon (x^2+y^2)^2
+\sqrt{2}(\sqrt{3}-1)(x^2+y^2)y
-3\sqrt{3}(\sqrt{3}-1)(x^2+y^2) \\
&\qquad 
+\sqrt{2}(\sqrt{3}-1)y
-\sqrt{3}+2=0, \\
n=4 &\colon (x^2+y^2)^2
+2(\sqrt{2}-1)(x^2+y^2)y
-6\sqrt{2}(\sqrt{2}-1)(x^2+y^2) \\
&\qquad 
+2(\sqrt{2}-1)y
-2\sqrt{2}+3=0, \\
n=5 &\colon 
(x^2+y^2)^2
+(\sqrt{2}(\sqrt{5}+1)^{\frac{1}{2}}5^{\frac{1}{4}}-\sqrt{5}-1)(x^2+y^2)x \\
&\qquad
+3(\tfrac{1}{\sqrt{2}}(\sqrt{5}+1)^{\frac{3}{2}}5^{\frac{1}{4}}-\sqrt{5}-5)(x^2+y^2) \\
&\qquad 
+(\sqrt{2}(\sqrt{5}+1)^{\frac{1}{2}}5^{\frac{1}{4}}-\sqrt{5}-1)x \\
&\qquad 
-\tfrac{1}{\sqrt{2}}(\sqrt{5}+1)^{\frac{3}{2}}5^{\frac{1}{4}}+\sqrt{5}+4=0.
\end{align*}
It is also easy to get the polar form, or view the curve as the intersection of the sphere with a parabolic cylinder
\begin{align*}
n=3 &\colon 2\sqrt{3}\xi_3^2+\sqrt{2}\xi_2+\xi_3-\sqrt{3}=0, \\
n=4 &\colon 2\sqrt{2}\xi_3^2+\xi_2+\xi_3-\sqrt{2}=0, \\
n=5 &\colon 2\sqrt{2}\cdot 5^{\frac{1}{4}}\xi_3^2+{\textstyle \sqrt{\sqrt{5}-1}}\xi_2+{\textstyle \sqrt{\sqrt{5}+1}}\xi_3-\sqrt{2}\cdot 5^{\frac{1}{4}}=0.
\end{align*}

\begin{figure}[htp]
\centering
\begin{tikzpicture}[>=latex,scale=5]

% tetrahedron

\draw[gray!40]
	({(sqrt(3)-1)/(sqrt(2))},0) -- (0,0) -- (0,{(sqrt(3)-1)/(sqrt(2))});
	
\draw[samples=100, smooth, gray!40]
	plot[domain=90:180] ({\x}:{( % A3
		sqrt(
			(sqrt(6*cos(\x)*cos(\x)+2)+sqrt(2)*cos(\x))*
			(sqrt(6*cos(\x)*cos(\x)+2)+sqrt(2)*cos(\x))+4
			) 
		- sqrt(6*cos(\x)*cos(\x)+2) - sqrt(2)*cos(\x)
			)/2})
	plot[domain=0:-90] ({\x}:{( % B3
		sqrt(
			(sqrt(6*sin(\x)*sin(\x)+2)+sqrt(2)*sin(\x))*
			(sqrt(6*sin(\x)*sin(\x)+2)+sqrt(2)*sin(\x))+4
			) 
		- sqrt(6*sin(\x)*sin(\x)+2) - sqrt(2)*sin(\x)
			)/2})
	plot[domain=180:270] ({\x}:{( % C3
		sqrt(
			(sqrt(2)*(2*cos(\x)*sin(\x)-1)/(cos(\x)+sin(\x)))*
			(sqrt(2)*(2*cos(\x)*sin(\x)-1)/(cos(\x)+sin(\x)))+4
			) 
		- (sqrt(2)*(2*cos(\x)*sin(\x)-1)/(cos(\x)+sin(\x)))
			)/2});

% a=b

\draw[gray!40]
	(-0.8,-0.8) -- (0.4,0.4);

\draw
	(-0.71,-0.71) -- (0,0);
	
% a=c
 
\foreach \a in {1,-1}
{
\begin{scope}[xscale=\a, rotate=45-45*\a]

\foreach \b in {1,-1}			
\draw[gray!40, yscale=\b, samples=100, smooth]
	plot[domain=0.1607:0.5176] (
		{ -(\x*\x*\x*\x-3*(3-sqrt(3))*\x*\x-sqrt(3)+2)/((sqrt(6)-sqrt(2))*(\x*\x+1))},
		{ -sqrt(\x*\x-((\x*\x*\x*\x-3*(3-sqrt(3))*\x*\x-sqrt(3)+2)*(\x*\x*\x*\x-3*(3-sqrt(3))*\x*\x-sqrt(3)+2))/((sqrt(6)-sqrt(2))*(sqrt(6)-sqrt(2))*(\x*\x+1)*(\x*\x+1))) }
		)
	(0.516,-0.028) -- ({(sqrt(3)-1)/(sqrt(2))},0)
	(-0.1606,0) -- (-0.16,0.013);	
			
\draw[samples=100, smooth]
	plot[domain=0.1607:0.5176] (
		{ -(\x*\x*\x*\x-3*(3-sqrt(3))*\x*\x-sqrt(3)+2)/((sqrt(6)-sqrt(2))*(\x*\x+1))},
		{ -sqrt(\x*\x-((\x*\x*\x*\x-3*(3-sqrt(3))*\x*\x-sqrt(3)+2)*(\x*\x*\x*\x-3*(3-sqrt(3))*\x*\x-sqrt(3)+2))/((sqrt(6)-sqrt(2))*(sqrt(6)-sqrt(2))*(\x*\x+1)*(\x*\x+1))) }
		)
	(0.516,-0.028) -- ({(sqrt(3)-1)/(sqrt(2))},0)
	(-0.1606,0) -- (-0.16,-0.013);

\draw[yscale=-1, samples=100, smooth]
	plot[domain=0.1607:0.268] (
		{ -(\x*\x*\x*\x-3*(3-sqrt(3))*\x*\x-sqrt(3)+2)/((sqrt(6)-sqrt(2))*(\x*\x+1))},
		{ -sqrt(\x*\x-((\x*\x*\x*\x-3*(3-sqrt(3))*\x*\x-sqrt(3)+2)*(\x*\x*\x*\x-3*(3-sqrt(3))*\x*\x-sqrt(3)+2))/((sqrt(6)-sqrt(2))*(sqrt(6)-sqrt(2))*(\x*\x+1)*(\x*\x+1))) }
		)
	(-0.1606,0) -- (-0.16,-0.013);	

\end{scope}
}
	
\node[rotate=45] at (-0.39,-0.45) {$a=b$};
\node at (0.2,-0.35) {$a=c$};	
\node[rotate=-60] at (-0.3,0) {$b=c$};

% octahedron

\begin{scope}[xshift=1.1cm]

\draw[gray!40]
	({sqrt(3)-sqrt(2)},0) -- (0,0) -- (0,{sqrt(2)-1});
	
\draw[samples=100, smooth, gray!40]
	plot[domain=90:180] ({\x}:{sqrt(cos(\x)*cos(\x)+2)-sqrt(cos(\x)*cos(\x)+1)}) % A4
	plot[domain=0:-90] ({\x}:{ % B4
		sqrt(
			(sqrt(2*sin(\x)*sin(\x)+2)+sin(\x))*
			(sqrt(2*sin(\x)*sin(\x)+2)+sin(\x))+1
			) 
		- sqrt(2*sin(\x)*sin(\x)+2) - sin(\x)
			})
	plot[domain=180:270] ({\x}:{ % C4
		sqrt(
			((cos(\x)*sin(\x)-sqrt(2))/(cos(\x)+sqrt(2)*sin(\x)))*
			((cos(\x)*sin(\x)-sqrt(2))/(cos(\x)+sqrt(2)*sin(\x)))+1
			) 
		- (cos(\x)*sin(\x)-sqrt(2))/(cos(\x)+sqrt(2)*sin(\x)) 
			});

% a=b

\draw[gray!40, shift={( {-sqrt(2)*(2+sqrt(3))}, {sqrt(3)*(2+sqrt(3))} )}]
	(-48.3:{2*sqrt(9+5*sqrt(3))}) arc (-48.3:-55.5:{2*sqrt(9+5*sqrt(3))}); 
	
\draw[shift={( {-sqrt(2)*(2+sqrt(3))}, {sqrt(3)*(2+sqrt(3))} )}]
	(-50.28:{2*sqrt(9+5*sqrt(3))}) arc (-50.28:-54.84:{2*sqrt(9+5*sqrt(3))}); 
	
% a=c

\foreach \a in {1,-1}			
\draw[gray!40, yscale=\a, samples=100, smooth] 
	plot[domain=0.105:0.3178] (
		{ -(\x*\x*\x*\x-6*(3-sqrt(6))*\x*\x-2*sqrt(6)+5)/(2*(sqrt(3)-sqrt(2))*(\x*\x+1))},
		{ -sqrt(\x*\x-((\x*\x*\x*\x-6*(3-sqrt(6))*\x*\x-2*sqrt(6)+5)*(\x*\x*\x*\x-6*(3-sqrt(6))*\x*\x-2*sqrt(6)+5))/((2*(sqrt(3)-sqrt(2))*(2*(sqrt(3)-sqrt(2)))*(\x*\x+1)*(\x*\x+1))) }
		)
	(0.314,-0.04) -- ({sqrt(3)-sqrt(2)},0)
	(-0.103,0) -- (-0.101,0.03);

\draw[samples=100, smooth] 
	plot[domain=0.105:0.3178] (
		{ -(\x*\x*\x*\x-6*(3-sqrt(6))*\x*\x-2*sqrt(6)+5)/(2*(sqrt(3)-sqrt(2))*(\x*\x+1))},
		{ -sqrt(\x*\x-((\x*\x*\x*\x-6*(3-sqrt(6))*\x*\x-2*sqrt(6)+5)*(\x*\x*\x*\x-6*(3-sqrt(6))*\x*\x-2*sqrt(6)+5))/((2*(sqrt(3)-sqrt(2))*(2*(sqrt(3)-sqrt(2)))*(\x*\x+1)*(\x*\x+1))) }
		)
	(0.314,-0.04) -- ({sqrt(3)-sqrt(2)},0)
	(-0.103,0) -- (-0.101,0.03);

\draw[yscale=-1, samples=100, smooth] 
	plot[domain=0.105:0.176] (
		{ -(\x*\x*\x*\x-6*(3-sqrt(6))*\x*\x-2*sqrt(6)+5)/(2*(sqrt(3)-sqrt(2))*(\x*\x+1))},
		{ -sqrt(\x*\x-((\x*\x*\x*\x-6*(3-sqrt(6))*\x*\x-2*sqrt(6)+5)*(\x*\x*\x*\x-6*(3-sqrt(6))*\x*\x-2*sqrt(6)+5))/((2*(sqrt(3)-sqrt(2))*(2*(sqrt(3)-sqrt(2)))*(\x*\x+1)*(\x*\x+1))) }
		)
	(-0.103,0) -- (-0.101,0.03);
	
% b=c

\foreach \a in {1,-1}			
\draw[xscale=\a, gray!40, samples=100, smooth] 
	plot[domain=0.133:0.4142] (
		{ -sqrt(\x*\x-((\x*\x*\x*\x-6*(2-sqrt(2))*\x*\x-2*sqrt(2)+3)*(\x*\x*\x*\x-6*(2-sqrt(2))*\x*\x-2*sqrt(2)+3))/((2*(sqrt(2)-1)*(2*(sqrt(2)-1))*(\x*\x+1)*(\x*\x+1))) },
		{ -(\x*\x*\x*\x-6*(2-sqrt(2))*\x*\x-2*sqrt(2)+3)/(2*(sqrt(2)-1)*(\x*\x+1))}
		)
	(0.02,0.413) -- (0,{sqrt(2)-sqrt(1)})
	(0,{2-sqrt(2)+sqrt(3)-sqrt(6)}) -- ++(-0.03,0.002);

\draw[samples=100, smooth] 
	plot[domain=0.133:0.4142] (
		{ -sqrt(\x*\x-((\x*\x*\x*\x-6*(2-sqrt(2))*\x*\x-2*sqrt(2)+3)*(\x*\x*\x*\x-6*(2-sqrt(2))*\x*\x-2*sqrt(2)+3))/((2*(sqrt(2)-1)*(2*(sqrt(2)-1))*(\x*\x+1)*(\x*\x+1))) },
		{ -(\x*\x*\x*\x-6*(2-sqrt(2))*\x*\x-2*sqrt(2)+3)/(2*(sqrt(2)-1)*(\x*\x+1))}
		)
	(-0.02,0.413) -- (0,{sqrt(2)-sqrt(1)})
	(0,{2-sqrt(2)+sqrt(3)-sqrt(6)}) -- ++(-0.03,0.002);

\draw[xscale=-1, samples=100, smooth] 
	plot[domain=0.133:0.222] (
		{ -sqrt(\x*\x-((\x*\x*\x*\x-6*(2-sqrt(2))*\x*\x-2*sqrt(2)+3)*(\x*\x*\x*\x-6*(2-sqrt(2))*\x*\x-2*sqrt(2)+3))/((2*(sqrt(2)-1)*(2*(sqrt(2)-1))*(\x*\x+1)*(\x*\x+1))) },
		{ -(\x*\x*\x*\x-6*(2-sqrt(2))*\x*\x-2*sqrt(2)+3)/(2*(sqrt(2)-1)*(\x*\x+1))}
		)
	(0,{2-sqrt(2)+sqrt(3)-sqrt(6)}) -- ++(-0.03,0.002);
		
\node[rotate=37] at (-0.18,-0.27) {$a=b$};
\node at (0.11,-0.23) {$a=c$};	
\node[rotate=-40] at (-0.2,-0.08) {$b=c$};	

\end{scope}

% icosahedron

\begin{scope}[xshift=1.8cm]

\draw[gray!40]
	({(sqrt(15)-sqrt(5)+sqrt(3)-3)/2},0) -- (0,0) -- (0,{( sqrt(10+2*sqrt(5)) -sqrt(5)-1 )/2});
	
\draw[samples=100, smooth, gray!40]
	plot[domain=90:180] ({\x}:{( % A5
		sqrt( 
			(sqrt((10-2*sqrt(5))*cos(\x)*cos(\x)+2*sqrt(5)+6)
				-(sqrt(5)-1)*cos(\x))*
			(sqrt((10-2*sqrt(5))*cos(\x)*cos(\x)+2*sqrt(5)+6)
				-(sqrt(5)-1)*cos(\x)) + 4 
			)  
		-(sqrt((10-2*sqrt(5))*cos(\x)*cos(\x)+2*sqrt(5)+6)
				-(sqrt(5)-1)*cos(\x)) 
			)/2})
	plot[domain=0:-90] ({\x}:{( % B5
		sqrt( 
			(sqrt((10+2*sqrt(5))*sin(\x)*sin(\x)+6*sqrt(5)+14)
				+(sqrt(5)+1)*sin(\x))*
			(sqrt((10+2*sqrt(5))*sin(\x)*sin(\x)+6*sqrt(5)+14)
				+(sqrt(5)+1)*sin(\x)) + 4 
			)  
		-(sqrt((10+2*sqrt(5))*sin(\x)*sin(\x)+6*sqrt(5)+14)
				+(sqrt(5)+1)*sin(\x)) 
			)/2})
	plot[domain=180:270] ({\x}:{ % C5
		sqrt(
			((sqrt(5)-1)*(cos(\x)*sin(\x)-sqrt(5)-2)/(2*cos(\x)+(sqrt(5)+1)*sin(\x)))*
			((sqrt(5)-1)*(cos(\x)*sin(\x)-sqrt(5)-2)/(2*cos(\x)+(sqrt(5)+1)*sin(\x))) + 1
			) 
		- ((sqrt(5)-1)*(cos(\x)*sin(\x)-sqrt(5)-2)/(2*cos(\x)+(sqrt(5)+1)*sin(\x)))
			});

\pgfmathsetmacro{\mm}{sqrt(195-6*sqrt(5))/58};

% a=b

\draw[gray!40, shift={( { (sqrt(5)-11)*\mm - sqrt(5) }, { (5*sqrt(5)+3)*\mm + 3 } )}]
	(-54:{ sqrt( 4*(13*sqrt(5)+2)*\mm +30 ) }) arc (-54:-59:{ sqrt( 4*(13*sqrt(5)+2)*\mm +30 ) }); 
	
\draw[shift={( { (sqrt(5)-11)*\mm - sqrt(5) }, { (5*sqrt(5)+3)*\mm + 3 } )}]
	(-55.13:{ sqrt( 4*(13*sqrt(5)+2)*\mm +30 ) }) arc (-55.13:-58.4:{ sqrt( 4*(13*sqrt(5)+2)*\mm +30 ) }); 

%% a=c

\foreach \a in {1,-1}			
\draw[yscale=\a, gray!40, samples=100, smooth] 
	plot[domain=0.061:0.1844] (
		{ -( \x*\x*\x*\x+3*(2*sqrt(15)-3*sqrt(5)+4*sqrt(3)-9)*\x*\x-2*sqrt(15)+3*sqrt(5)-4*sqrt(3)+8 )/((sqrt(5)-sqrt(3))*((sqrt(3)-sqrt(1))*(\x*\x+1))},
		{ -sqrt(\x*\x-((\x*\x*\x*\x+3*(2*sqrt(15)-3*sqrt(5)+4*sqrt(3)-9)*\x*\x-2*sqrt(15)+3*sqrt(5)-4*sqrt(3)+8)*(\x*\x*\x*\x+3*(2*sqrt(15)-3*sqrt(5)+4*sqrt(3)-9)*\x*\x-2*sqrt(15)+3*sqrt(5)-4*sqrt(3)+8)/((sqrt(5)-sqrt(3))*((sqrt(3)-sqrt(1))*(sqrt(5)-sqrt(3))*((sqrt(3)-sqrt(1))*(\x*\x+1)*(\x*\x+1))) }
		)
	(0.1844,0) -- ++(-0.003,0.03)
	(-0.061,0) -- ++(0.0013,0.015);

\draw[samples=100, smooth] 
	plot[domain=0.061:0.1844] (
		{ -( \x*\x*\x*\x+3*(2*sqrt(15)-3*sqrt(5)+4*sqrt(3)-9)*\x*\x-2*sqrt(15)+3*sqrt(5)-4*sqrt(3)+8 )/((sqrt(5)-sqrt(3))*((sqrt(3)-sqrt(1))*(\x*\x+1))},
		{ -sqrt(\x*\x-((\x*\x*\x*\x+3*(2*sqrt(15)-3*sqrt(5)+4*sqrt(3)-9)*\x*\x-2*sqrt(15)+3*sqrt(5)-4*sqrt(3)+8)*(\x*\x*\x*\x+3*(2*sqrt(15)-3*sqrt(5)+4*sqrt(3)-9)*\x*\x-2*sqrt(15)+3*sqrt(5)-4*sqrt(3)+8)/((sqrt(5)-sqrt(3))*((sqrt(3)-sqrt(1))*(sqrt(5)-sqrt(3))*((sqrt(3)-sqrt(1))*(\x*\x+1)*(\x*\x+1))) }
		)
	(0.1844,0) -- ++(-0.003,-0.03)
	(-0.061,0) -- ++(0.0013,-0.015);
			
\draw[yscale=-1, samples=100, smooth] 
	plot[domain=0.061:0.105] (
		{ -( \x*\x*\x*\x+3*(2*sqrt(15)-3*sqrt(5)+4*sqrt(3)-9)*\x*\x-2*sqrt(15)+3*sqrt(5)-4*sqrt(3)+8 )/((sqrt(5)-sqrt(3))*((sqrt(3)-sqrt(1))*(\x*\x+1))},
		{ -sqrt(\x*\x-((\x*\x*\x*\x+3*(2*sqrt(15)-3*sqrt(5)+4*sqrt(3)-9)*\x*\x-2*sqrt(15)+3*sqrt(5)-4*sqrt(3)+8)*(\x*\x*\x*\x+3*(2*sqrt(15)-3*sqrt(5)+4*sqrt(3)-9)*\x*\x-2*sqrt(15)+3*sqrt(5)-4*sqrt(3)+8)/((sqrt(5)-sqrt(3))*((sqrt(3)-sqrt(1))*(sqrt(5)-sqrt(3))*((sqrt(3)-sqrt(1))*(\x*\x+1)*(\x*\x+1))) }
		)
	(-0.061,0) -- ++(0.0013,-0.015);
		
%% b=c

\foreach \a in {1,-1}			
\draw[xscale=\a, gray!40, samples=100, smooth]
	plot[domain=0.093:0.284] (
		{ -sqrt( \x*\x -( ( 100*\x*\x*\x*\x-324.2102*\x*\x+8.0701 )*( 100*\x*\x*\x*\x-324.2102*\x*\x+8.0701 ) )/( 3228.0361*(\x*\x+1)*(\x*\x+1) ) 
				) },
		{ -( 100*\x*\x*\x*\x-324.2102*\x*\x+8.0701 )/( 56.8158*(\x*\x+1) ) }
		)
	(0,0.2841) -- ++(0.027,-0.002)
	(0,-0.0925) -- ++(0.015,+0.001);
	 
\draw[samples=100, smooth] 
	plot[domain=0.093:0.19] (
		{ -sqrt( \x*\x -( ( 100*\x*\x*\x*\x-324.2102*\x*\x+8.0701 )*( 100*\x*\x*\x*\x-324.2102*\x*\x+8.0701 ) )/( 3228.0361*(\x*\x+1)*(\x*\x+1) ) 
				) },
		{ -( 100*\x*\x*\x*\x-324.2102*\x*\x+8.0701 )/( 56.8158*(\x*\x+1) ) }
		)
	(0,-0.0925) -- ++(0.015,+0.001);

\draw[xscale=-1, samples=100, smooth]
	plot[domain=0.093:0.158] (
		{ -sqrt( \x*\x -( ( 100*\x*\x*\x*\x-324.2102*\x*\x+8.0701 )*( 100*\x*\x*\x*\x-324.2102*\x*\x+8.0701 ) )/( 3228.0361*(\x*\x+1)*(\x*\x+1) ) 
				) },
		{ -( 100*\x*\x*\x*\x-324.2102*\x*\x+8.0701 )/( 56.8158*(\x*\x+1) ) }
		)
	(0,-0.0925) -- ++(0.015,+0.001);
			
\node[rotate=32] at (-0.1,-0.18) {\footnotesize $a=b$};
\node at (0.09,-0.15) {\footnotesize $a=c$};	
\node[rotate=-35] at (-0.13,-0.06) {\footnotesize $b=c$};	

\end{scope}

\node at (-0.2,-1.2) 
	{\includegraphics[scale=0.16]{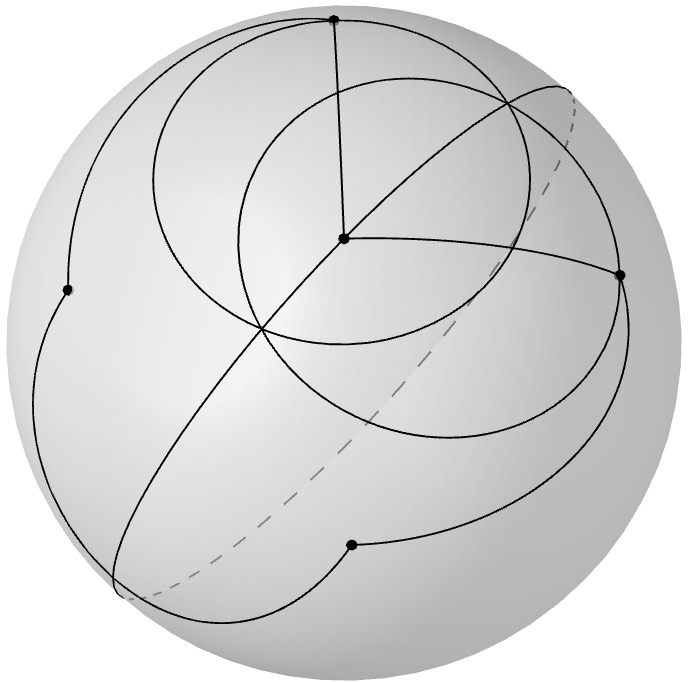}};

\node at (0.7,-1.2) 
	{\includegraphics[scale=0.20]{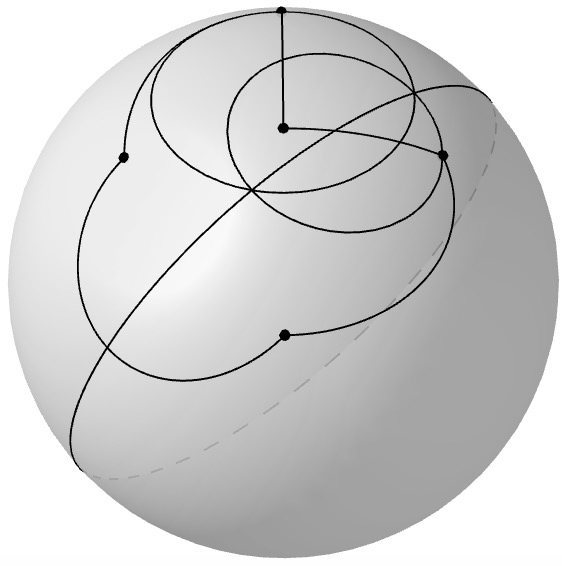}};
	
\node at (1.6,-1.2) 
	{\includegraphics[scale=0.16]{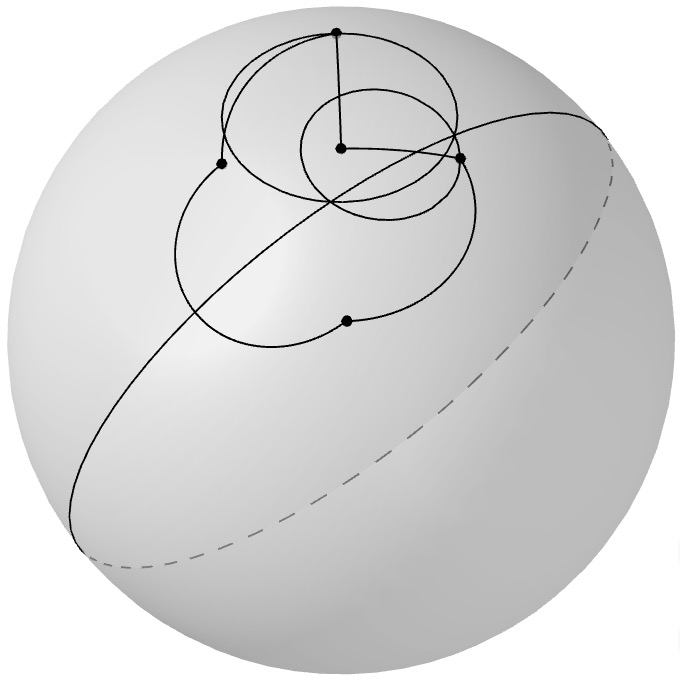}};
		
\end{tikzpicture}
\caption{Reduction curves in the moduli space.}
\label{reduce}
\end{figure}

Figure \ref{reduce} is the perspective picture of the reduced curves in the $M$-projection and in the sphere. From the picture, we notice that the curve $b=c$ is tangential to and is also below the boundary $\gamma_A$ of the moduli space for $n=4$. We confirmed this through Taylor expansion. The picture also suggests that each pentagonal subdivision has a unique equilateral pentagon (i.e., $a=b=c$). The uniqueness of the equilateral pentagonal subdivisions are rigorously verified in \cite{ay1,awy}.

%%% reduction

\subsection{Size of Moduli Space}

In a stereographic projection, we consider the fan region between the origin and a curve $r=r(\theta)$, $\theta\in[\alpha,\beta]$. The corresponding region in the sphere has area
\[
\int_{\alpha}^{\beta}\frac{2r^2d\theta}{1+r^2}.
\]
Then the area of the region given by \eqref{eqE} and \eqref{eqD} is ($\theta=\phi-\alpha$)
\[
\int_{\alpha}^{\frac{1}{2}\pi}
\frac{2r^2d\theta}{1+r^2}
=\int_{\alpha}^{\frac{1}{2}\pi} 
\left(1-\frac{R}{\sqrt{R^2+1}}\right)d\theta
=\tfrac{1}{2}\pi-\alpha
-F(\tfrac{1}{2}\pi-\alpha,\lambda^{-1}i),
\]
where $R=\lambda\sec\theta$, and
\[
F(t,k)
=\int_0^{t}
\frac{dt}{\sqrt{1-k^2\sin^2t}}
\]
is the elliptic integral of the first kind.

We may apply the formula to the $A$-projection of $\gamma_A$ to get the area of the moduli space part in $\Omega_5$
\[
A_5=\tfrac{1}{6}\pi
-F(\tfrac{1}{6}\pi,\lambda_A^{-1}i).
\]
We may also apply the formula to the $B$-projection of $\gamma_B$ to get the area of the moduli space part in $\Omega_{13}$ (for $n=5$, the notation $A_{13}$ also includes the moduli space part in $\Omega_{19}$)
\[
A_{13}=\left(\tfrac{1}{2}-\tfrac{1}{n}\right)\pi
-F((\tfrac{1}{2}-\tfrac{1}{n})\pi,\lambda_B^{-1}i).
\]

We may use the same idea to get (for $n=5$, the notation $A_8$ also includes the moduli space part in $\Omega_{14}$)
\begin{align*}
A_4
&=\tfrac{1}{6}\pi
-F(\tfrac{1}{6}\pi,\lambda_C^{-1}i), \\
A_8
&=\left(\tfrac{1}{2}-\tfrac{1}{n}\right)\pi
-F((\tfrac{1}{2}-\tfrac{1}{n})\pi,\lambda_C^{-1}i).
\end{align*}
We may also use the fact that $\gamma_C$ has the same formula in $\Omega_4$ and $\Omega_8$ to get
\[
A_2+A_4+A_8
=\tfrac{1}{2}\pi-F(\tfrac{1}{2}\pi,\lambda_C^{-1}i).
\]

Finally, the moduli space contains the whole triangles $\Omega_1,\Omega_2,\Omega_3,\Omega_7$, which are congruent triangles and have the same area
\[
A_1=A_2=A_3=A_7=\frac{1}{6N}4\pi,\quad
N=4,8,20\text{ for }n=3,4,5.
\]
The whole area of the moduli space is 
\[
A_1+A_2+A_3+A_4+A_5+A_7+A_8+A_{13}
=\begin{cases}
0.8600517493 \pi, & n=3, \\
0.4602931496 \pi, & n=4, \\
0.1954959087 \pi, & n=5.
\end{cases}
\]

\end{document}